\documentclass[oneside,english,reqno]{amsart}
\usepackage{palatino}
\usepackage[T1]{fontenc}
\usepackage[latin1]{inputenc}
\usepackage{amssymb}
\IfFileExists{url.sty}{\usepackage{url}}
                      {\newcommand{\url}{\texttt}}

\makeatletter

\providecommand{\LyX}{L\kern-.1667em\lower.25em\hbox{Y}\kern-.125emX\@}
\newcommand{\noun}[1]{\textsc{#1}}

 \theoremstyle{plain}    
 \newtheorem*{conjecture*}{Conjecture} 
 \theoremstyle{definition}
 \newtheorem*{defn*}{Definition}
 \theoremstyle{plain}    
 \newtheorem{lem}{Lemma} 
 \theoremstyle{plain}    
 \newtheorem*{lem*}{Lemma} 
 \theoremstyle{remark}
 \newtheorem*{rem*}{Remark}
 \theoremstyle{plain}    
 \newtheorem{thm}{Theorem} 

\DeclareMathOperator{\LE}{LE}

\DeclareMathOperator{\cov}{cov}

\newcommand{\brmul}{\discretionary{\mbox{$\,\cdot$}}{}{}}

\usepackage{babel}
\makeatother
\begin{document}

\title{loop-erased random walk on a torus in dimensions 4 and above}

\author{Itai Benjamini}

\author{Gady Kozma}

\maketitle

\section{Introduction}

A well known phenomenon in probabilistic constructions in $\mathbb{R}^{d}$
or $\mathbb{Z}^{d}$ is that usually some critical dimension $d$
exists, above which the geometry of $\mathbb{R}^{d}$ ceases to play
any significant role, and the process behaves like a similar non-geometric
object, such as a tree, a complete graph, etc. Usually, this also
corresponds to {}``mean field behavior'', a term meaning that for
the random variables of interest one has $\mathbb{E}X^{n}\approx (\mathbb{E}X)^{n}$.
At the critical dimension itself, mean field behavior is also expected,
but when compared to the non-geometric object one gets a logarithmic
correction.

Many results confirming this general philosophy exist. See \cite{HS90}
for results about percolation, \cite{HS92} for results about the
self-avoiding walk, \cite{DS98} for results about lattice trees,
and \cite{S95} for a general survey. In particular, the problem of
loop-erased random walk on $\mathbb{Z}^{d}$ is well studied.

Loop-erased random walk is a process that starts from a random walk
on some graph and then removes all loops in chronological order, or
in other words, whenever the random walk hits the partial path, the
loop just created is erased and the process continues. The result
is a random simple path. Originally \cite{L80} suggested as a model
for the self-avoiding walk (a random walk conditioned not to hit itself),
better understanding of its structure has situated it as an important
object in combinatorics and mathematical physics. See \cite{S00}
for a survey, and the complementary \cite{LSW}. For a survey with
a different focus, see \cite{L99}. Note also the recent \cite{BKPS}
--- the uniform spanning tree is an object closely related to loop-erased
random walk, but the structure of its phase transitions in various
dimensions is richer. For other recent results of interest, see \cite{BLPS01,K,LPS}.

It is well known that the critical dimension of loop-erased random
walk on $\mathbb{Z}^{d}$ is $4$, since above this dimension a random
walk does not intersect itself enough and the process of loop-erasure
is local and uninteresting. See \cite[chapter 7]{L96}. Further, loop-erased
walk is one of the few models where the logarithmic correction is
known precisely, with a correction of $\log ^{-1/3}$ loop-erased
random walk on $\mathbb{Z}^{4}$ is similar to the regular random
walk on $\mathbb{Z}^{4}$, see \cite{L95}.

With so much known, it seems strange that a small change in settings
could provoke significant difficulties. To understand why, let us
examine the question we are interested in precisely. Let $T$ be a
discrete torus, $\mathbb{Z}^{d}/(N\mathbb{Z})^{d}$ for some large
$N$. Let $b$ and $e$ be two points on a torus, and let $R$ be
a random walk starting from $b$ and stopped on $e$. We wish to say
something about the loop-erasure of $R$. The results for $\mathbb{Z}^{d}$
all use the fact that the random walk does not intersect itself enough.
However, in our settings the random walk does a very long walk ---
of the order of $N^{d}$ --- in a relatively small space, and intersects
itself over and over again. Thus it is definitely \textbf{not} true
that the random walk and its loop erasure are similar! The random
walk is essentially a random set that covers a large portion of the
torus. Its loop-erasure is much thinner --- as we will see, the expected
size is $N^{d/2}$.

The geometry-less model we have in mind is the complete graph. There
are a number of ways this model can be analyzed, but our favorite
is using the notion of the Laplacian random walk. A Laplacian random
walk from $b$ to $e$, two points on an arbitrary graph $G$, is
constructed inductively by solving, at each step, the discrete Dirichlet
problem \begin{equation}
f(e)=1,\; f|_{\gamma }\equiv 0,\; \Delta f|_{G\setminus (\gamma \cup \{e\})}\equiv 0\label{eq:laplac}\end{equation}
where $\gamma $ is the partially constructed path and $\Delta $
is the discrete Laplacian. The walk then continues to the next point
using $f$ as weights. This model was suggested in \cite{LEP86} and
was shown to be equivalent to loop-erased random walk in \cite{L87}.
The case of the complete graph is very easy to analyze, since if the
partially constructed curve $\gamma $ has length $i$ then \[
f(v)=\begin{cases}
 0 & v\in \gamma \\
 1 & v=e\\
 \frac{1}{i+1} & \textrm{otherwise}\end{cases}.\]
and then the probability of the walk to terminate in the next step
is $\frac{i}{N}$. This gives a closed formula\[
\mathbb{P}(\#\LE (R)=k)=\frac{k-1}{N}\prod _{i=1}^{k-2}1-\frac{i}{N}\quad .\]
In particular, we see that the correct scaling is $\sqrt{N}$ and
that $\#\LE (R)/\sqrt{N}$ converges to a limiting distribution with
density $te^{-t^{2}/2}$. Unfortunately, we do not know how to analyze
more interesting graphs using the Laplacian random walk, nor can we
show the existence of a limiting distribution for $\#\LE (R)$ on,
say, the torus.

Thus we have a good basis to claim that mean field behavior in our
case should be $\sqrt{|T|}=N^{d/2}$. For $d<4$ this does not happen
--- indeed known results for $d=1$ (trivial) and $d=2$ (\cite{K00,K00b},
see also \cite{LSW}) and computer simulations for $d=3$ \cite{GB90}
show that even a single branch of the loop-erased walk is too big%
\footnote{We believe that the growth exponents in $d=2,3$ are the same on $\mathbb{Z}^{d}$
and $T_{N}^{d}$, but this is beyond the scope of this paper.%
}. We shall show that mean field behavior does occur for $d>4$. In
the critical dimension itself, we can only show an upper bound, and
we do not calculate the precise logarithmic correction (we do have
some good evidence for a conjecture on the precise logarithmic correction
needed --- $\log ^{1/6}N$ --- see page \pageref{con:log16}). Namely,
our results are

\begin{thm}
\label{thm:upper}If $d>4$ then a loop-erased random walk $L$ on
the $(N,d)$-torus starting from a point $b$ and stopped when hitting
a point $e$ has the estimate\[
\mathbb{P}(\#L>\lambda N^{d/2})\leq Ce^{-c\lambda }\quad .\]
If $d=4$ then \[
\mathbb{P}(\#L>\lambda N^{2+\epsilon })\leq Ce^{-c\lambda }\quad \forall \epsilon >0.\]
Where the constants $C$ and $c$ may depend on $d$ and on $\epsilon $.
\end{thm}

\begin{thm}
\label{thm:dge5}Let $d\geq 5$. Let $b$ be a point in $T=T_{N}^{d}$
and let $e$ be a random, uniform point in $T$. Let $R$ be a random
walk on $T$ starting from $b$ and stopped at $e$. Let $\lambda \geq N^{-1/2}$.
Then \[
\mathbb{P}(\#\LE (R)\leq \lambda N^{d/2})\leq C\lambda \log \lambda ^{-1}\]

\end{thm}
Returning to the cases of $d\leq 3$, we see that the reason for non-mean-field
behavior is strong local intersections and these \emph{increase} the
size of the loop-erased walk. Therefore we are tempted to conjecture

\begin{conjecture*}
Let $G$ be a vertex transitive finite graph, and let $b$ and $e$
be two random points in $G$. Let $R$ be a random walk starting from
$b$ and stopped when hitting $e$. Then\[
\mathbb{E}\#\LE (R)\geq c\sqrt{|G|}\]

\end{conjecture*}
A graph $G$ is vertex transitive when, for every two vertices $v$
and $w$ there exists a graph automorphism of $G$ carrying $v$ to
$w$. The requirement that $G$ is vertex transitive is supported
by the standard {}``extreme non-transitive'' example of a tree of
size $N$, where the loop-erased random walk between $b$ and $e$
is of course the \emph{only} path between $b$ and $e$ and its length
is bounded by $C\log N$.

We wish to end this introduction with one last conjecture. Returning
to the analysis of the complete graph using the Laplacian random walk,
we note that this analysis does not change by much if one considers
$\alpha $-power Laplacian random walk, which is a walk one gets if
one takes as weights for any step the function $f^{\alpha }$ where
$f$ is defined by (\ref{eq:laplac}) --- this generalization was
also discussed in \cite{LEP86}. For the complete graph we get that
the size of a typical path is $N^{1/(1+\alpha )}$. We ask: is this
behavior replicated in a $d$-dimensional torus for $d>d_{\alpha }^{\textrm{crit}}$?

\begin{conjecture*}
For $\alpha \leq 1$, $d_{\alpha }^{\textrm{crit}}=\frac{2(1+\alpha )}{\alpha }$,
i.e.~for any $d>\frac{2(1+\alpha )}{\alpha }$ the typical path of
a $\alpha $-weighted Laplacian random walk on a $d$-dimensional
torus is of size $N^{d/(1+\alpha )}$ while for smaller $d$'s this
does not hold.
\end{conjecture*}
We have no good conjecture on the value of the critical dimension
for $\alpha >1$, though it does seem (again, we have no proof of
that) that for $\alpha =\infty $ (which corresponds to a non-probabilistic
process which simply proceeds to the point where $f$ attains its
maximum) the process gives a straight line from $b$ to $e$ in all
dimensions, so one might say the critical dimension is $1$.

We wish to thank Chris Hoffman, Dan Romik and Oded Schramm for useful
discussions.

\subsection{About the proof}

The basic question behind the solution is {}``what is the probability
of a random walk of length $L$ will hit a loop-erased walk of length
$L$?'' (in dimension $4$ we need to differentiate between these
two lengths, but only by a sub-polynomial factor). When the probability
is larger then some constant $c>0$, then this is the $L$ we seek,
as this means that the probability of a loop-erased random walk to
go further than $\lambda L$ is exponentially small in $\lambda $.
Since a loop-erased walk is a complicated object, let us first ask
{}``what is the probability of a random walk of length $L$ will
hit some set $\Omega $ of size $L$?'' This probability is largest
when $\Omega $ is rather spread out. Take as an example $\Omega $
to be a random collection of points on the torus. It is easy to calculate
the expected number of intersections of a random walk with $\Omega $
and the second moment and to derive from both the estimate that the
probability is $\approx L^{2}N^{-d}$ so this gives that the $L$
we look for is $<CN^{d/2}$. When the set $\Omega $ is rather dense
--- for example, for a ball --- a similar calculation will give that
$L<CN^{d/2+1}$. The difference between $\Omega $ dense and sparse
manifests itself in the calculation of the second moment --- see (\ref{eq:expx2_1})-(\ref{eq:expx2_2})
below.

However, $\Omega $ is not a general set but a loop-erased random
walk. The arguments we sketched above can be done locally, and we'd
get that in every ball of radius $r$, the loop-erased random walk
is not much larger than $r^{d/2+1}$. Effectively, this means that
the set is spread out. We take this estimate and plug it directly
into the calculation of the second moment and get a much better estimate
for the intersection probability. Thus the proof is \emph{recursive},
getting better estimates at each step. For $d>4$ two or at most three
steps are necessary to get the true estimate, $N^{d/2}$. This argument
is done in lemma \ref{lem:main}.

\subsection{Reading recommendations}

Section \ref{sec:upper} is probably the one deserving most attention.
While the main ideas are sketched above, the devil is in the details
and the interested reader might want to read through the proof and
do the {}``exercise'' --- not so designated explicitly --- of simplifying
the proof with a cost of $\sqrt{\log }$ in the final result. Section
\ref{sec:Absolute-times} is technical and most readers would probably
agree that the conclusion (theorem \ref{thm:absolute}) is not surprising.
The proof of lemma \ref{lem:abs} is the core --- as for lemma \ref{lem:times},
you might opt to read its statement but skip its proof. And again,
verify that the claim is trivial if one is willing to lose a factor
of $\log $ (the argument is contained in the first half-page of the
proof of lemma \ref{lem:abs}). Section \ref{sec:The-lower-bound}
contains the proof of theorem \ref{thm:dge5} and is quite short.
While there are alternative, more complicated approaches that might
prove a little more we have not included them. There are some comments
and hints at the end of section \ref{sec:The-lower-bound} --- we
hope they make at least some sense. We have collected some well known
and unsurprising facts we use (and their proofs) in the appendix.
We hope this makes the paper more accessible to non-experts and students.
Lemmas with numbers like {}``A.7'' are to be found in the appendix.

\subsection{Standard notations }

In the sequel we denote by $C$ and $c$ positive constants which
may depend on the dimension but on nothing else. $C$ will usually
pertain to constants which are {}``large enough'' and $c$ to constants
which are {}``small enough''. The notation $x\approx y$ is a short
hand for $cx\leq y\leq Cx$. In dimension $4$ we shall prove only
imprecise estimates, namely that the length of the loop-erased walk
is $<N^{2+\epsilon }$. All constants $C$ and $c$ may depend on
this $\epsilon $ as well. Similarly, all constants implicit in notations
such as $O$ and $\approx $ might depend on $d$ and $\epsilon $.
Occasionally we shall number constants for clarity. When we write
$\log x$ we always mean $\max \{\log x,1\}$ and $\log 0=1$.

The $(N,d)$-torus, denoted by $T_{N}^{d}$ is the set $\mathbb{Z}^{d}/(N\mathbb{Z})^{d}$
endowed with the graph structure derived from $\mathbb{Z}^{d}$ and
the distance derived from the $l_{2}$ norm on $\mathbb{Z}^{d}$.
The distance of $v$ and $w$ will be denoted by $|v-w|$ while distance
of sets will be denoted by $d(\cdot ,\cdot )$. A ball of radius $r$
and center $v$ in either $\mathbb{Z}^{d}$ or $T_{N}^{d}$ will be
denoted by $B(v,r)$ and its inner boundary (namely, all points in
$B$ with an edge leading outside of $B$) by $\partial B(v,r)$.

\section{\label{sec:upper}The upper bound}

We will need to examine the effect of adding a section to a path and
how it might increase the length of its loop-erasure. We shall always
assume that the section we add starts at $0$, so that we are looking
at a path $\gamma :\{-m,\dotsc ,n\}\rightarrow T$ and define, in
addition to the usual loop-erasure of $\gamma $, which we will denote
by $\LE (\gamma )$, the \textbf{continued} loop-erasure, which we
shall denote by $\LE ^{+}(\gamma )$. Here are both definitions:

\begin{defn*}
For a finite path $\gamma :\{-m,\dotsc ,n\}\rightarrow T$ in a graph
$T$ we define its loop erasure, $\LE (\gamma )$, which is a simple
path in $T$, by the consecutive removal of loops from $\gamma $.
Formally, \begin{eqnarray*}
\LE (\gamma )_{0} & := & \gamma (-m)\\
\LE (\gamma )_{i+1} & := & \gamma (j_{i}+1)\quad j_{i}:=\max \{j:\gamma (j)=\LE (\gamma )_{i}\}\quad .
\end{eqnarray*}
Naturally, this is defined for all $i$ such that $j_{i}<n$. The
continued loop-erasure is a subset of $\LE (\gamma )$ defined by\begin{alignat*}{2}
\LE ^{+}(\gamma )_{i} & :=\LE (\gamma )_{I+i}\quad  & I & :=\min \{i:j_{i}\geq 0\}
\end{alignat*}
The notations $\LE (\gamma [A,B])$ and $\LE ^{+}(\gamma [A,B])$
stand for the loop-erasure and continued loop-erasure of the segment
of $\gamma $ going from $A$ to $B$. When we write $-\infty $ in
place of $A$ we just mean the beginning of the path, nothing more.
\end{defn*}

\begin{defn*}
Let $d\geq 4$ be the dimension and let $N\in \mathbb{N}$. Let $R$
be a path in $T=T_{N}^{d}$ such that the negative part is fixed and
the positive part is a random walk on $T$. Let $b=R(0)$. Let $v\in T$
and $0<r<\frac{1}{8}N$, and assume for simplicity that $b\not \in B(v,2r)$.
Let $t_{i}$ be stopping times defined by $t_{0}=0$ and then inductively
\begin{eqnarray}
t_{2i+1} & := & \min \{t\geq t_{2i}:R(t)\in \partial B(v,2r)\}\label{eq:defti}\\
t_{2i} & := & \min \{t\geq t_{2i-1}:R(t)\in \partial B(v,4r)\}\quad .\nonumber 
\end{eqnarray}
Let $f:\mathbb{R}\rightarrow \mathbb{R}$ be an increasing function.
Then we say that the ($d$-dimensional) random walk has the $f$\textbf{-property}
if one has\begin{equation}
\mathbb{P}(\#(\LE ^{+}(R[-\infty ,t_{i}])\cap B(v,r))>\lambda f(r)\, |\, R[t_{2j},t_{2j+1}]\forall j)\leq Ce^{-c\lambda }\label{eq:property}\end{equation}
which should hold for every such $v$ and $r$, every $\lambda >0$,
every $i\in \mathbb{N}$ and any path we put in the negative portion
of $R$.

The conditioning here, in words, is on any arbitrary set of paths
between $t_{2j}$ and $t_{2j+1}$, and in particular on the points
$R(t_{i})$ themselves. Notice that we do not condition on the value
of the $t_{i}$'s.
\end{defn*}
Let us remark that for the proof of the upper bound it is enough to
consider the case where $R$ has no negative part, and then $\LE ^{+}\equiv \LE $. 

\begin{lem}
\label{lem:main}Let $d\geq 4$. Then
\begin{enumerate}
\item If the $d$-dimensional random walk satisfies the $r^{\alpha }\log ^{\beta }r$-property
for $r^{\alpha }\brmul \log ^{\beta }r\gg r^{d-2}\log ^{-3}r$ then
it also satisfies the $r^{\alpha /2+1}\log ^{(\beta +3)/2}r$-property.
\item If the $d$-dimensional random walk satisfies the $r^{d-2}\log ^{-3}r$-property
then it also satisfies the $r^{d/2}\sqrt{\log \log r}$-property.
\item If it satisfies the $r^{\alpha }\log ^{\beta }$-property for $r^{\alpha }\log ^{\beta }r\ll r^{d-2}\log ^{-3}r$
then it satisfies the $r^{d/2}$-property.
\end{enumerate}
\end{lem}
Case 2 is not really necessary for the proof of the theorem, we include
it here mainly for completeness.

\begin{proof}
Denote the function given to us (e.g.~$r^{\alpha }\log ^{\beta }r$)
by $f(r)$ and the result (e.g.\ $r^{\alpha /2+1}\log ^{(\beta +3)/2}r$)
by $g(r)$. Let $t_{i}$ be the stopping times from the definition
of the $f$-property. The main part of the lemma will consider the
events in $R[t_{i},t_{i+1}]$ for some particular odd $i$. Therefore
let us fix $i>0$. Denote $L_{i,v,r}:=\#(\LE ^{+}(R[0,t_{i}])\cap B(v,r))$.
Clearly $L_{2i+1,v,r}\leq L_{2i,v,r}$ so if we prove the lemma for
all $i$ odd it will also hold for $i$ even. To fix notations, we
consider the time span $[-\infty ,t_{i}]$ as the {}``past'' and
$\left]t_{i},t_{i+1}\right]$ is the {}``present''. 

We start by examining the past. Let $w\in B(v,r)$ and $s\leq \frac{1}{8}r$.
The first step is to show that (\ref{eq:property}) holds if we replace
the ball but keep the stopping times, i.e\begin{equation}
\mathbb{P}(\#(\LE ^{+}(R[0,t_{i}])\cap B(w,s))>\lambda f(s)\, |\, R[t_{2j},t_{2j+1}]\forall j)\leq Ce^{-c\lambda }\quad .\label{eq:inner}\end{equation}
We generalize the notation $L_{i,v,r}$ to $L_{i,w,s}:=\#(\LE ^{+}(R[0,t_{i}])\cap B(w,s))$,
that is, again, the loop-erased random walk inside a smaller ball
measured at the stopping times pertaining to the larger ball.

Here our conditioning by everything outside the ball is crucial. Let
$K_{j}\in \mathbb{N}$ be some arbitrary numbers, and let $\gamma _{j,k}$
be paths ($1\leq k\leq K_{j}$) in $B(v,4r)\setminus B(w,2s)$ such
that $\gamma _{j,1}$ is a path going from $R(t_{2j-1})\in \partial B(v,2r)$
to $\partial B(w,2s)$, $\gamma _{j,k}$ for $1<k<K_{j}$ is a path
from $\partial B(w,4s)$ to $\partial B(w,2s)$ and $\gamma _{j,K_{j}}$
is a path from $\partial B(w,4s)$ to $R(t_{2j})\in \partial B(v,4r)$.
If $K_{j}=1$ then let $\gamma _{j,1}$ be a path from $R(t_{2j-1})$
to $R(t_{2j})$. Then we can sum over all such combinations of $K$
and $\gamma $ as follows. Denote by $X$ the event $L_{i,w,s}>\lambda f(s)$.
Let $Y_{K,\gamma }$ be the event that for all $j$, the random walk
on $[t_{2j-1},t_{2j}]$ follows $\gamma _{j,1}$ until $\partial B(w,2s)$,
then stays within $B(w,4s)$, then follows $\gamma _{j,2}$ etc. until
finally exiting from $B(v,4r)$. Then \begin{align}
\mathbb{P}(X\, |\, R[t_{2j},t_{2j+1}]\forall j) & =\sum _{K,\gamma }\mathbb{P}(X\, |\, R[t_{2j},t_{2j+1}]\forall j\cap Y_{K,\gamma })\cdot \mathbb{P}(Y_{K,\gamma }\, |\, R[t_{2j},t_{2j+1}]\forall j)\nonumber \\
 & \leq Ce^{-c\lambda }\sum _{K,\gamma }\mathbb{P}(Y_{K,\gamma }\, |\, R[t_{2j},t_{2j+1}]\forall j)=Ce^{-c\lambda }\quad .\label{eq:BvrBws}
\end{align}
Of course, we used the $f$-property for $w$, $s$ and the index
$\sum _{j=1}^{(i-1)/2}K_{j}$; and the fact that $B(w,4s)\subset B(v,2r)$.

The inequality (\ref{eq:inner}) is not useful as it should be since
most balls of radius $s$ (for $s\ll r$) are empty anyway. However,
another consequence of the conditioning is the fact that (\ref{eq:inner})
is independent from the event $\LE ^{+}(R[-\infty ,t_{i}])\cap B(w,4s)=\emptyset $.
The reason is that $L_{i,w,4s}=0$ if and only if the segment inside
$B(w,4s)$ is cut {}``from the root'', i.e.~for some $u_{1}<u_{2}<\dotsb <u_{2n}$,
$n\in \{1,2,\dotsc \}$ we must have $R[u_{2i-1},u_{2i}]\cap B(w,4s)=\emptyset $
and $R(u_{2i})=R(u_{2i+1})$. Whether this happens in the positive
or negative part of $R$ is immaterial --- in both cases this is an
event that happens outside $B(w,4s)$ therefore it is an event we
condition on. We get \begin{equation}
\mathbb{P}(L_{i,w,s}>\lambda f(s)\, |\, R[t_{2j},t_{2j+1}]\forall j)\leq Ce^{-c\lambda }\mathbb{P}(L_{i,w,4s}\neq 0\, |\, R[t_{2j},t_{2j+1}]\forall j)\quad .\label{eq:bignz}\end{equation}

Let $\gamma =\gamma _{i}$ be the (chronologically) first $G$ elements
of $\LE ^{+}(R[-\infty ,t_{i}])\cap B(v,r)$ where $G$ is some number.
If $\LE ^{+}(R[-\infty ,t_{i}])\cap B(v,r)$ contains less than $G$
elements, take $\gamma =\LE ^{+}(R[-\infty ,t_{i}])\cap B(v,r)$.
(\ref{eq:inner}) and (\ref{eq:bignz}) allow us to get a {}``second-order
estimate for $\gamma $''. By this we mean the quantity\[
V_{s}:=\#\{w_{1},w_{2}\in \gamma :|w_{1}-w_{2}|\leq s\}\]
which has the estimate\begin{equation}
\mathbb{P}(\#\gamma >\delta \mathbb{E}L_{i,v,r}\textrm{ and }V_{s}>\lambda \log (s/\delta )f(s)\#\gamma )\leq Ce^{-c\lambda }\label{eq:Vs2}\end{equation}
for any parameters $\lambda >0$ and $0<\delta <1$.

Before starting the proof of (\ref{eq:Vs2}) let us just remark that
the first condition and the variable $\delta $ are unfortunate technicalities.
The {}``essentials'' of (\ref{eq:Vs2}) are really the stronger
claim $\mathbb{P}(V_{s}>\lambda (\log s)f(s)\#\gamma )\leq Ce^{-c\lambda }$,
but we don't know how to prove it. Also note that it is rather easy
to show $\mathbb{P}(V_{s}>\lambda (\log r)f(s)\#\gamma )\leq Ce^{-c\lambda }$,
saving us all the mucking with $\delta $ later on, but this inequality
will cost us a $\sqrt{\log r}$ in the final result of theorem \ref{thm:upper}.
\begin{proof}
[Proof of (\ref{eq:Vs2})]Cover $B(v,r)$ by balls $\{B_{j}\}$ of
radius $2s$ such that any two points of distance $\leq s$ are inside
at least one $B_{j}$, and such that each point is covered $\leq C$
times. Examine one $B_{j}=B(w_{j},2s)$. We have (not writing the
{}``$|\, R[t_{2j},t_{2j+1}]\forall j$'' for brevity) \[
\mathbb{E}L_{i,v,r}>c\sum _{j}\mathbb{P}(L_{i,w_{j},8s}>0)\stackrel{(\ref {eq:bignz})}{\geq }ce^{c\lambda }\sum _{j}\mathbb{P}(L_{i,w_{j},2s}>\lambda f(s))\quad \forall \lambda .\]
Denote by $X_{\mu }$ the total volume of the balls $B_{j}$ where
$L_{i,w_{j},2s}>\mu f(s)$ and get $\mathbb{E}X_{\mu }\leq Ce^{-c\mu }s^{d}\mathbb{E}L_{i,v,r}$.
This gives, using $\mathbb{P}(X_{\mu }>e^{c\mu }\mathbb{E}X_{\lambda })\leq e^{-c\mu }$,
\[
\mathbb{P}(X_{\mu }>Cs^{d}e^{-c\mu }\mathbb{E}L_{i,v,r})\leq e^{-c\mu }\quad \forall \mu \]
and shoving in $\#\gamma $ in a way that might look, for now, a little
artificial, we get\[
\mathbb{P}(\#\gamma >\delta \mathbb{E}L_{i,v,r}\textrm{ and }X_{\mu }>Cs^{d}\delta ^{-1}e^{-c\mu }\#\gamma )\leq e^{-c\mu }\quad \forall \mu .\]
Taking $\mu _{k}=\lambda \log (s/\delta )+Ck$ and assuming that $\lambda >C$
for some $C$ sufficiently large (as we may, without loss of generality),
we get \begin{equation}
\mathbb{P}(\#\gamma >\delta \mathbb{E}L_{i,v,r}\textrm{ and for some }k,\textrm{ }X_{\mu _{k}}>c^{-k}\#\gamma )\leq Ce^{-c\lambda }\quad .\label{eq:xmusmall}\end{equation}
Now, since $V_{s}\leq \sum _{j}\#(\gamma \cap B_{j})\cdot L_{i,w_{j},2s}$,
then\[
V_{s}\leq \#\gamma (\lambda \log (s/\delta )+C)f(s)+\sum _{k=1}^{\infty }X_{\mu _{k}}(\lambda \log (s/\delta )+Ck)f(s).\]
If it happens that $X_{\mu _{k}}\leq c^{-k}\#\gamma $ for all $k$
i.e.~the opposite of the second half of the event in (\ref{eq:xmusmall}),
then \begin{align*}
V_{s} & \leq \lambda \log (s/\delta )f(s)\#\gamma +\sum _{k=1}^{\infty }(c^{-k}\#\gamma )f(s)(\lambda \log (s/\delta )+Ck)\\
 & \leq C\lambda \log (s/\delta )f(s)\#\gamma 
\end{align*}
and we get (\ref{eq:Vs2}). This argument works for any $s\leq \frac{1}{16}r$
but (\ref{eq:Vs2}) holds for larger $s$ too (there's not much point
in $s>2r$ of course) --- we only have to pay in the constant $C$.
\end{proof}
We want (\ref{eq:Vs2}) to hold not for one particular $s$ but for
all $s$ and the simplest version of such an inequality is\begin{equation}
\mathbb{P}(\#\gamma >\delta \mathbb{E}L_{i,v,r}\textrm{ and }\exists s\textrm{ s.t. }V_{s}>\lambda \log ^{2}(s/\delta )f(s)\#\gamma )\leq Ce^{-c_{1}\lambda }\label{eq:Vs3}\end{equation}
which follows from using (\ref{eq:Vs2}) with $\lambda _{s}:=\lambda \log (s/\delta )$
and summing over $s$. 

Continuing the proof of the lemma, it is now time to examine the present.
We keep the notations of $G$, $\gamma $ and $V_{s}$. For an odd
$i$ we want to estimate the probability \[
p_{i}:=\mathbb{P}(R[t_{i},t_{i+1}]\cap \gamma \neq \emptyset )\quad .\]
Lemma \ref{lem:cond} allows us to consider a unconditioned random
walk starting from $R(t_{i})$ and stopped on $\partial B(v,4r)$
instead of $R$. Denote it by $R'$. Denote by $X_{i}$ the number
of intersections of $R'$ with $\gamma $, so $p_{i}\approx \mathbb{P}(X_{i}>0)$.
We have\[
\mathbb{E}(X_{i}\, |\, \textrm{past})=\sum _{t=t_{i}}^{t_{i+1}}\sum _{w\in \gamma }\mathbb{P}(R'(t)=w\, |\, \textrm{past})\quad .\]
For $r^{2}\leq t-t_{i}\leq 2r^{2}$ we have for half of the $w\in B(v,r)$
that $\mathbb{P}(\{R'(t)=w\}\cap \{t<t_{i+1}\})>cr^{-d}$ ({}``half
of the $w$'s'' means that we need $t-t_{i}+||w-R(t_{i})||_{1}$
to be even, otherwise the probability is zero). Therefore \begin{equation}
\mathbb{E}(X_{i}\, |\, \textrm{past})>cr^{2-d}\#\gamma \quad .\label{eq:expx}\end{equation}

Next estimate $\mathbb{E}(X_{i}^{2}\, |\, \textrm{past})$. Assume
until further notice that $V_{s}\leq \lambda \log ^{2}(s/\delta )\brmul f(s)\#\gamma $
for some $\delta $ and $\lambda $ and for all $s$. Then

\begin{align}
\mathbb{E}(X_{i}^{2}\, |\, \textrm{past}) & =\sum _{t_{1},t_{2},w_{1},w_{2}}\mathbb{P}(R'(t_{i})=w_{i})\leq \label{eq:expx2_1}\\
 & \leq 2\sum _{\Delta =0}^{\infty }\sum _{k=0}^{\infty }\sum _{t,w_{1},w_{2}}\mathbb{P}\Big (R(t)=w_{1},R(t+\Delta )=w_{2},\nonumber \\
 & \qquad \qquad \qquad \qquad \qquad k\sqrt{\Delta }\leq |w_{1}-w_{2}|<(k+1)\sqrt{\Delta }\Big )\nonumber 
\end{align}
Examine one couple of $w_{1},w_{2}\in \gamma $ with $k\sqrt{\Delta }\leq |w_{1}-w_{2}|$.
Remembering the independence of the past from the present we can estimate
the probability of one summand with a standard estimate on the end
point of a random walk of length $\Delta $ starting from $w_{1}$.
We get\[
\mathbb{P}(R(t)=w_{1},R(t+\Delta )=w_{2})\leq Cr^{-d}\Delta ^{-d/2}e^{-k^{2}/2}\quad .\]
We sum over all $t$. Since, easily, $\mathbb{P}(t_{i+1}-t_{i}>\Delta )\leq Ce^{-c\Delta /r^{2}}$
and since $\mathbb{E}(t_{i+1}-t_{i}\, |\, t_{i+1}-t_{i}>\Delta )\leq C\max \{r^{2},\Delta \}$
we get\[
\sum _{t}\mathbb{P}(R(t)=w_{1},R(t+\Delta )=w_{2})\leq Ce^{-c\Delta /r^{2}}\max \{r^{2},\Delta \}r^{-d}\Delta ^{-d/2}e^{-k^{2}/2}\quad .\]
Plugging this into (\ref{eq:expx2_1}) we get\begin{equation}
\mathbb{E}(X_{i}^{2}\, |\, \textrm{past})\leq C\sum _{\Delta =0}^{\infty }\sum _{k=0}^{\infty }e^{-c\Delta /r^{2}}\max \{r^{2},\Delta \}r^{-d}\Delta ^{-d/2}e^{-k^{2}/2}V_{(k+1)\sqrt{\Delta }}\label{eq:expx2_2}\end{equation}
For all our functions $f$ (that is, all the specific functions we
named in the statement of the lemma) we have \begin{align*}
\sum _{k=0}^{\infty }e^{-k^{2}/2}V_{(k+1)\sqrt{\Delta }} & \leq \lambda \#\gamma \sum _{k=1}^{\infty }e^{-(k-1)^{2}/2}f(k\sqrt{\Delta })\log ^{2}(k\sqrt{\Delta }/\delta )\leq \\
 & \leq C\lambda \#\gamma f(\sqrt{\Delta })\log ^{2}(\Delta /\delta )\quad .
\end{align*}
Similarly, for all our functions $f$ we have\begin{eqnarray}
\lefteqn{\sum _{\Delta =0}^{\infty }e^{-c\Delta /r^{2}}\max \{r^{2},\Delta \}\Delta ^{-d/2}f(\sqrt{\Delta })\log ^{2}(\Delta /\delta )\leq } &  & \label{eq:sumf}\\
 & \qquad \qquad \qquad \qquad  & \leq Cr^{2}\sum _{\Delta =0}^{r^{2}}\Delta ^{-d/2}f(\sqrt{\Delta })\log ^{2}(\Delta /\delta )\quad .\nonumber 
\end{eqnarray}
(\ref{eq:expx2_2}) and (\ref{eq:sumf}) give\[
\mathbb{E}(X_{i}^{2}\, |\, \textrm{past})\leq C\lambda r^{2-d}\#\gamma \sum _{\Delta =0}^{r^{2}}\Delta ^{-d/2}f(\sqrt{\Delta })\log ^{2}(\Delta /\delta )\]
and then with (\ref{eq:expx}) and the standard inequality $\mathbb{P}(X>0)\geq (\mathbb{E}X)^{2}/\mathbb{E}X^{2}$
we get\begin{equation}
\mathbb{P}(X_{i}>0\, |\, \textrm{past})>c\frac{r^{2-d}\#\gamma }{\lambda \sum \Delta ^{-d/2}f(\sqrt{\Delta })\log ^{2}(\Delta /\delta )}\quad .\label{eq:crux}\end{equation}
This inequality is the heart of the proof. We recall that we assumed
$V_{s}\leq \lambda \brmul \log ^{2}(s/\delta )\brmul f(s)\#\gamma $
to get it.

Fix $G=\mu g(r)$ where $\mu >1$ is some variable which we will fix
later and where $g$ is as defined in the beginning of the lemma.
Let $H=2\left\lfloor g(r)r^{-2}\right\rfloor $ where $\left\lfloor \cdot \right\rfloor $
is the integer value. Let $X_{1}=X_{1}(\mu )$ be the event that $\#\gamma _{i}=G$,
let $X_{2}=X_{2}(\lambda ,\delta ,\mu )$ be the event that $V_{s}\leq \lambda \log ^{2}(s/\delta )f(s)G$
for all $s$ ($\lambda $ and $\delta $ are two additional variables)
and let $X_{3}=X_{3}(\mu )$ be the event that $R[t_{j+1}-t_{j}]\cap \gamma =\emptyset $
for all odd $i\leq j\leq i+H$. The events comprising $X_{3}$ are
(conditioning on the $R(t_{j})$) independent, therefore we may use
(\ref{eq:crux}) $\frac{1}{2}H$ times to get\begin{equation}
\mathbb{P}(X_{3}\, |\, X_{1}\cap X_{2})\leq \left(1-c\frac{r^{2-d}\mu g(r)}{\lambda \sum _{\Delta =1}^{r^{2}}\Delta ^{-d/2}f(\sqrt{\Delta })\log ^{2}(\Delta /\delta )}\right)^{\frac{1}{2}H}\leq 1-\frac{c\mu }{\lambda \log ^{2}\delta ^{-1}}\quad .\label{eq:whyg}\end{equation}
To see the rightmost inequality in (\ref{eq:whyg}), for each of the
cases in the formulation of the lemma, apply the corresponding $f$
and $g$ and estimate the sum. Indeed, (\ref{eq:whyg}) is the inequality
that governs the connection between $f$ and $g$. Note that the formulation
of the lemma is a little lax: if $f(r)=r^{\alpha }\log ^{\beta }r$
with $\alpha >d-2$ then we can actually prove the lemma with $g=r^{\alpha /2}\log ^{(\beta +2)/2}$
i.e.~one $\sqrt{\log r}$ factor better than the formulation of the
lemma. This additional $\sqrt{\log r}$ factor is here only for the
case $\alpha =d-2$ and $\beta <-3$. Have no fear --- this factor
will disappear in the conclusion of theorem \ref{thm:upper}.

The proof of the lemma will now follow by induction over $i$. We
use a {}``jumping induction'' that assumes that for some $k$ and
$K$ we have the inequality $\mathbb{P}(L_{i,v,r}>\nu g(r))\leq Ke^{-k\nu }$
for all $\nu >0$ and then proves the same for $L_{i+H,v,r}$ (the
case $i=0$ needs no explanation). Therefore we need first to calculate
how much $L_{i,v,r}$ can change in between. Clearly, if $R([t_{j},t_{j+1}])$
does not intersect $\LE ([R[0,t_{j}])$ then \[
L_{j+1,v,r}-L_{j,v,r}\leq t_{j+1}-t_{j}\quad .\]
These variables have the simple estimate\begin{equation}
\mathbb{P}(t_{j+1}-t_{j}>\nu r^{2})\leq Ce^{-c\nu }\label{eq:tijump}\end{equation}
irrespectively of $R(t_{j+1})$ and $R(t_{j})$ for all $j$ odd.
Denote by $A_{i}$ the sum of $\frac{1}{2}H$ of those, and get a
similar estimate (see lemma \ref{lem:sumexpexp}):\begin{equation}
\mathbb{P}\left(A_{i}>\nu g(r)\right)\leq Ce^{-c_{2}\nu }\quad A_{i}:=\sum _{\substack{ j=i\\
 j\textrm{ odd}}
}^{i+H}t_{j+1}-t_{j}\quad .\label{eq:Ai}\end{equation}

Next we make the following important assumption:\begin{equation}
G>\delta \mathbb{E}L_{i,v,r}\quad \forall i.\label{eq:assump}\end{equation}
Actually, we want it to be true independently of the value of $\mu $,
so we really need $g(r)>\delta \mathbb{E}L_{i,v,r}$. This holds for
$\delta $ sufficiently small, but it is inconvenient to fix the value
of $\delta $ at this point, as it depends on some constants (depending
on $d$ only) which are determined only later. Therefore we shall
perform the necessary calculations with $\delta $ a variable and
finally fix its value as some constant when we have all the information
at hand, see (\ref{eq:espval}). With a value of $\delta $ satisfying
(\ref{eq:espval}), or smaller, (\ref{eq:assump}) will hold. 

It is time to compare $L_{i,v,r}$ with $L_{i+H,v,r}$. $L_{i+H,v,r}$
might be larger than $\nu g(r)$ for the simple reason that $A_{i}$
is very large. Let $\tau \leq \nu $ be yet another variable describing
what {}``very large'' means and we may estimate this phenomenon
simply by \[
\sum _{n=\tau }^{\nu }\mathbb{P}(\{L_{i,v,r}>(\nu -n-1)g(r)\}\cap \{A_{i}>ng(r)\})\quad .\]
{}``Simply'' because we ignore any effect of intersections. If,
however, $A_{i}$ is not as large we need both $L_{i,v,r}$ to be
rather large, and $X_{3}$, i.e.~to have no intersections with a
path of length $G=\mu g(r)$ during the last $H$ {}``moves''. We
need to assume $\mu +\tau <\nu $ for this to make sense, and this
assumption holds until (\ref{eq:endmutaunu}) below and we will not
repeat it. All in all we get\begin{align*}
\mathbb{P}(L_{i+H,v,r} & >\nu g(r))\leq \mathbb{P}(\{L_{i,v,r}>(\nu -\tau )g(r)\}\cap X_{3})+\\
 & +\sum _{n=\tau }^{\nu }\mathbb{P}(\{L_{i,v,r}>(\nu -n-1)g(r)\}\cap \{A_{i}>ng(r)\})\quad \forall i,\nu ,\tau ,\mu 
\end{align*}
(the parameter $\mu $ hides in the definition of $X_{3}$). For the
first summand we have by (\ref{eq:Vs3}), (\ref{eq:assump}), (\ref{eq:whyg})
and the induction hypothesis that \begin{align*}
\lefteqn{\mathbb{P}(\{L_{i,v,r}>(\nu -\tau )g(r)\}\cap X_{3})\leq } & \\
 & \qquad \leq \mathbb{P}(\{L_{i,v,r}>(\nu -\tau )g(r)\}\setminus X_{2})+\mathbb{P}(\{L_{i,v,r}>(\nu -\tau )g(r)\}\cap X_{3}\cap X_{2})\\
 & \qquad \leq Ce^{-c_{1}\lambda }+Ke^{-k(\nu -\tau )}\left(1-\frac{c\mu }{\lambda \log ^{2}\delta ^{-1}}\right)\quad \forall i,\nu ,\tau ,\lambda ,\mu ,\delta 
\end{align*}
and estimating the other summands using (\ref{eq:Ai}) we get\begin{align}
\mathbb{P}(L_{i+H,v,r}>\nu g(r)) & \leq Ke^{-k(\nu -\tau )}\left(1-\frac{c\mu }{\lambda \log ^{2}\delta ^{-1}}\right)+Ce^{-c_{1}\lambda }+\nonumber \\
 & +\sum _{n=\tau }^{\nu }Ke^{-k(\nu -n-1)}\cdot Ce^{-c_{2}n}\quad \forall i,\nu ,\tau ,\lambda ,\mu ,\delta .\label{eq:endmutaunu}
\end{align}
Having arrived at this closed formula, we only need to pick our variables
carefully. First pick $\tau =\left\lfloor C\log \delta ^{-1}\right\rfloor $
for some $C$ sufficiently large. This will give, if $k<c_{2}/2$,
that\[
\sum _{n=\tau }^{\nu }Ke^{-k(\nu -n-1)}\cdot Ce^{-c_{2}n}\leq C\frac{e^{-C\log \delta ^{-1}}}{1-e^{-c_{2}/2}}Ke^{-k(\nu -\tau )}\leq C\delta Ke^{-k(\nu -\tau )}\quad .\]
Next we pick $\lambda =C\nu $ and $\mu =\frac{1}{2}\nu $, and the
requirement $\mu +\tau <\nu $ translates to $\nu >C\log \delta ^{-1}$.
We get from everything that\begin{align*}
\mathbb{P}(L_{i+H,v,r}>\nu g(r)) & \leq Ke^{-kv}\left(e^{kC\log \delta ^{-1}}\left(1-\frac{c}{\log ^{2}\delta ^{-1}}+C\delta \right)+Ce^{-c\nu }\right)
\end{align*}
Pick $k=c\log ^{-3}\delta ^{-1}$ and get, for $\delta $ sufficiently
small and $\nu >C\log \delta ^{-1}$ that \[
\mathbb{P}(L_{i+H,v,r}>\nu g(r))\leq Ke^{-k\nu }\left(1-\frac{c}{\log ^{2}\delta ^{-1}}\right)\quad .\]
Pick $K$ sufficiently large so that the inequality $\mathbb{P}(L_{i,v,r}>\nu g(r))\leq Ke^{-k\nu }$
will hold trivially for $\nu \leq C\log \delta ^{-1}$ --- notice
that because $k=c\log ^{-3}\delta ^{-1}$ we have that $K$ does not
depend on $\delta $ --- and our induction is complete. With these
$k$ and $K$, the inequality $\mathbb{P}(L_{i,v,r}>\nu g(r))\leq Ke^{-k\nu }$
is preserved from $i$ to $i+H$ and since it clearly holds for $i\leq H$
then it holds for all $i$.

Is this the end of the lemma? Almost. We still need to justify the
assumption (\ref{eq:assump}). The estimate \noun{$\mathbb{P}(L_{i,v,r}>\nu g(r))\leq Ke^{-k\nu }$}
gives $\mathbb{E}L_{i,v,r}\leq g(r)\frac{K}{k}\leq Cg(r)\log ^{3}\delta ^{-1}$.
Therefore (remember that $G>g(r)$) the assumption reduces to the
inequality\begin{equation}
g(r)>g(r)\cdot (C\delta \log ^{3}\delta ^{-1})\quad .\label{eq:espval}\end{equation}
Taking $\delta $ sufficiently small this will hold, and the lemma
is proved.
\end{proof}
\begin{lem}
\label{lem:mainpp}The $d$-dimensional random walk has the $f$-property
for \begin{equation}
f_{d}(r):=\begin{cases}
 r^{d/2} & d>4\\
 r^{2+\epsilon } & d=4\end{cases}\quad .\label{eq:deff}\end{equation}

\end{lem}
\begin{proof}
Trivially, the $d$-dimensional random walk has the $r^{d}$-property.
Therefore we may apply lemma \ref{lem:main} twice for $d>6$, thrice
for $d=6$ or $5$ and $\log \epsilon ^{-1}$ times for $d=4$.
\end{proof}

\begin{proof}
[Proof of theorem \ref{thm:upper}]Lemma \ref{lem:mainpp} gives\[
\mathbb{P}(L_{i,v,r}>\lambda f(r))\leq Ce^{-c\lambda }\]
where $L_{i,v,r}=\#(\LE (R[0,t_{i}])\cap B(v,r))$ for any $v$ and
$r$ satisfying $b\not \in B(v,2r)$, where $f$ is defined by (\ref{eq:deff}).
Note that at this point we do not need the formulation in terms of
continued process, and we may set the negative part of $R$ to empty.
If in addition $e\not \in B(v,4r)$ then the event that $R$ is stopped
between $t_{I}$ and $t_{I+1}$ is external to the ball, therefore
we get that (\ref{eq:deff}) holds for $I$. Since the section of
the walk from $t_{I}$ until the time when $R$ hits $e$ can only
decrease $\LE (R)\cap B(v,r)$ we get \[
\mathbb{P}(\#(L\cap B(v,r))>\lambda f(r))\leq Ce^{-c\lambda }\]
However, we can cover our torus by balls $B(v_{i,j},N2^{-i})$ with
the property $b,e\not \in B(v_{i,j},4N2^{-i})$ and with the number
of $j$'s corresponding to each $i$ bounded by a constant. Therefore
for some constant $c_{3}$ sufficiently small we have\begin{align}
\mathbb{P}(\#L>\lambda f(N)) & \leq \mathbb{P}(\exists i,j\textrm{ s.t. }L\cap B(v_{i,j},r)>c_{3}\lambda 2^{i/4}f(r))\label{eq:sumballs}\\
 & \leq \sum _{i=0}^{c\log N}Ce^{-c\lambda 2^{i/4}}\leq Ce^{-c\lambda }\quad .\qedhere \nonumber 
\end{align}

\end{proof}
\begin{rem*}
The same techniques can be improved to show that \[
\mathbb{P}(\#L>\lambda f(N))\leq Ce^{-c\lambda ^{2}}\]
where $f$ is given by (\ref{eq:deff}). The basic phenomenon behind
this estimate is that to get a path of length $\lambda f(N)$, we
need to have that each of the $\lambda $ sections of the random walk,
which are essentially independent, would not intersect any other.
Since there are $c\lambda ^{2}$ couples, the true estimate of the
probability is square-exponential, as above. The analysis required
to get this estimate is not inherently more difficult than that of
the exponential estimate, but is more technical and we decided to
represent the simpler exponential estimate.

On the other hand, we are not aware of a simpler version of the proof
that gives an estimate of the decay of the probability worse than
exponential. This follows from the recursive character of the proof.
Thus, lemma \ref{lem:main} may be simplified by removing the requirement
that the probability decays exponentially, but it then cannot be used
recursively to get a reasonable final result. Similarly, the very
strong independence condition in lemma \ref{lem:main}, that the probability
estimate inside every ball is independent of everything that happens
outside the ball, cannot be relaxed without destroying the ability
of the lemma to be used recursively.

We wish to reiterate that the only major simplification we are aware
of of this proof is the one discussed after (\ref{eq:Vs2}) (page
\pageref{eq:Vs2}). It saves the discussion after (\ref{eq:BvrBws}),
i.e.~the one leading to (\ref{eq:bignz}), as well as each and every
appearance of the parameter $\delta $. The cost is an added $\sqrt{\log }$
factor in the formulation of the theorem.
\end{rem*}
\begin{conjecture*}
\label{con:log16}The accurate upper bound in dimension $4$ is\[
N^{2}\log ^{1/6}N\quad .\]

\end{conjecture*}
The method above may be refined in many points and an estimate of
the type $N^{2}\log ^{\alpha }N$ may be achieved for rather small
$\alpha $'s. However, a fundamental difficulty is the fact that the
sum in the denominator of (\ref{eq:crux}) truly depends on $N$,
which means that the second moment methods used here alone cannot
give a precise result.

\section{\label{sec:Absolute-times}Absolute times}

The proof of the lower bound is, as will be seen in section \ref{sec:The-lower-bound},
quite simple once a good estimate of the upper bound is available.
Actually, one might think about the recursive nature of the proof
of the upper bound in the following terms: {}``the proof of the upper
bound was only possible once a good estimate of the upper bound was
available''.

Unfortunately, we were not able to get a reasonable proof of the lower
bound using only lemma \ref{lem:main}. The problem is that we need
to know what happens at absolute times, i.e.~to fix some $t$ and
get an estimate for $\LE (R[0,t])$. Calculations true for $t_{i}$
do not hold automatically for a fixed $t$. Apriori, one cannot rule
out behavior such as {}``the loop-erased random walk is much denser
if $t$ is divisible by $1024$'', since the $t_{i}$'s might avoid
those {}``bad absolute times''. The purpose of this section is to
show that this ridiculous behavior does not occur. 

The first step is to learn something about the distribution of the
$t_{i}$'s. Since $t_{i}$ is a sum of the return times to some sphere,
and these return times are more-or-less independent, we would expect
a central limit theorem. We don't need something so precise --- we
shall prove below (lemma \ref{lem:times}) a large deviation estimate
of the sort one would expect from a Gaussian variable, and this will
be enough. We start with

\begin{lem}
\label{lem:pseudo}Let $X_{1},\dotsc ,X_{n}$ be variables with the
properties\begin{align}
\mathbb{P}(|X_{i}|>\lambda \, | & \, X_{1},\dotsc ,X_{i-1},X_{i+1},\dotsc ,X_{n})\leq Ce^{-c\lambda }\label{eq:Xiexpo}\\
\mathbb{E}(X_{i_{1}}\dotsm X_{i_{k}}| & \, X_{i_{k+1}},\dotsc ,X_{i_{l}})\leq \prod _{j=1}^{k}C\exp (-c\min _{\substack{ 1\leq m\leq l\\
 m\neq j}
}|i_{j}-i_{m}|)\label{eq:Xipseudo}
\end{align}
where (\ref{eq:Xipseudo}) needs to hold only for $i_{1},\dotsc ,i_{l}$
all different. Then for all $\lambda <cn^{1/4}$\[
\mathbb{P}\left(\left|\sum X_{i}\right|>\lambda \sqrt{n}\right)\leq Ce^{-c\lambda ^{2}}\quad .\]

\end{lem}
We interpret the condition (\ref{eq:Xipseudo}) in the case $k=l=1$
as saying $\mathbb{E}X_{i}=0$ for all $i$. In the case $k>1$, we
call (\ref{eq:Xipseudo}) a {}``pseudo independence'' relation,
because, rather than claiming that $\mathbb{E}\prod X_{i}=0$, as
we would have for independent variables, we get that it is exponentially
small in the distance, so that if the $i_{k}$'s are relatively sparse,
it will be extremely small. Actually, it is possible to replace $\exp (-ck)$
with any sequence $a_{k}$ with $\sum a_{k}<C$.

The proof is a pretty standard exercise: a calculation (which can
be done either directly or by comparing to the case of independent
exponential variables) can show that for $k<c\sqrt{n}$,\[
\mathbb{E}\left(\sum X_{i}\right)^{2k}\leq (Ckn)^{k}\quad .\]
Taking $k=c\lambda ^{2}$ and using Markov's inequality will give
the lemma. We skip the gory details.

\begin{lem}
\label{lem:times}Let $b\in T_{N}^{d}$ and let $R$ be a random walk
on $T$ starting from $b$. Let $C<r<\frac{1}{8}N$, $v\in T_{N}^{d}$
and let $t_{i}$ be the stopping times defined by (\ref{eq:defti}).
Then there exists numbers $E=E(r)\approx N^{d}r^{2-d}$ and $\sigma =\sigma (r)\approx E$
such that\begin{equation}
\mathbb{P}(|t_{n}-nE|>\lambda \sigma \sqrt{n})\leq Ce^{-c\lambda ^{2}}\label{eq:explamb2}\end{equation}
for all $n\in \mathbb{N}$ and $\lambda <cn^{1/4}$.
\end{lem}
\begin{proof}
The point is of course to show that the variables $t_{i+1}-t_{i}$
are pseudo independent and apply lemma \ref{lem:pseudo}. The first
thing to note is that the distributions of $R(t_{i})$ converge exponentially.
Let $q_{1}$ and $q_{2}$ be two distributions on $\partial B(v,2r)$,
and denote\[
\epsilon :=\sum _{x\in B(v,2r)}|q_{1}(x)-q_{2}(x)|\quad .\]
 Let $R_{\mu }$, $\mu =1,2$ be random walks starting from a point
on $\partial B(v,2r)$ chosen with the distribution $q_{\mu }$ and
stopped when hitting $\partial B(v,4r)$. Let $p_{\mu }$ be the distributions
on the hit points of $R_{\mu }$. Then\begin{equation}
p_{1}(w)-p_{2}(w)=\sum _{x\in \partial B(v,2r)}(q_{1}(x)-q_{2}(x))\pi (x,w)\label{eq:p12q12pi}\end{equation}
where $\pi (x,w)$ is the probability of a random walk starting from
$x$ to hit $w$. Let $A^{+}\subset \partial B(v,2r)$ be the set
where $q_{1}(x)\geq q_{2}(x)$, and define\[
D^{+}(w)=\sum _{x\in A^{+}}|q_{1}(x)-q_{2}(x)|\pi (x,w)\quad .\]
Clearly\[
\sum _{w\in \partial B(v,4r)}D^{+}(w)=\sum _{x\in A^{+}}\sum _{w}|q_{1}(x)-q_{2}(x)|\pi (x,w)=\frac{1}{2}\epsilon \]
and similarly for $D^{-}$ defined equivalently using $A^{-}:=\partial B(v,2r)\setminus A^{+}$.
Furthermore, the inequality $\pi (x,w)\approx r^{1-d}$ (see lemma
\ref{lem:hitp}) gives that $D^{\pm }(w)\approx \epsilon r^{1-d}$
and therefore \[
|D^{+}(w)-D^{-}(w)|\leq (1-c)(D^{+}(w)+D^{-}(w))\]
for some constant $c>0$. This gives\begin{align}
\sum _{w\in \partial B(v,4r)}|p_{1}(w)-p_{2}(w)| & =\sum _{w}|D^{+}(w)-D^{-}(w)|\label{eq:expo}\\
 & \leq (1-c)\sum _{w}D^{+}(w)+D^{-}(w)=(1-c)\epsilon \nonumber 
\end{align}
and we see that the $L^{1}$ distance between the distributions has
contracted. An identical calculation works when the random walk starts
from $\partial B(v,4r)$ and stops at $\partial B(v,2r)$ (see the
remark following lemma \ref{lem:hitp}) therefore we see that there
is only one limiting distribution as $i$ increases, and that the
$L^{1}$ distance to this distribution decreases exponentially with
$i$. In other words, if $t_{i}^{\mu }$ are stopping times defined
by (\ref{eq:defti}) for the walks $R_{\mu }$ then we get\begin{equation}
\sum _{w}|\mathbb{P}(R_{1}(t_{i}^{1})=w)-\mathbb{P}(R_{2}(t_{i}^{2})=w)|\leq \epsilon e^{-ci}\quad .\label{eq:l1expo}\end{equation}
This $L^{1}$ estimate allows to get a uniform estimate for every
$w$ and $i>0$:\begin{equation}
|\mathbb{P}(R_{1}(t_{i}^{1})=w)-\mathbb{P}(R_{2}(t_{i}^{2})=w)|\leq Ce^{-ci}\min _{\mu =1,2}\mathbb{P}(R_{\mu }(t_{i}^{\mu })=w).\label{eq:linfexpo}\end{equation}
Indeed, take the distributions of $R_{i-1}^{\mu }$ as the $q_{\mu }$'s
in (\ref{eq:p12q12pi}) and together with (\ref{eq:l1expo}) and $\pi (x,w)\leq Cr^{1-d}$
get that \[
|\mathbb{P}(R_{1}(t_{i}^{1})=w)-\mathbb{P}(R_{2}(t_{i}^{2})=w)|\leq Cr^{1-d}e^{-ci}\quad .\]
 In the other direction, $\pi (x,w)\geq cr^{1-d}$ gives $\mathbb{P}(R_{\mu }(t_{i}^{\mu })=w)\geq cr^{1-d}$
and we get (\ref{eq:linfexpo}).

To make notations simpler, let $B_{i}$ be $\partial B(v,2r)$ if
$i$ is odd and $\partial B(v,4r)$ if $i$ is even. Now, each $t_{i+1}-t_{i}$
has an exponential distribution%
\footnote{For $i$ even, $t_{i+1}-t_{i}$ has a rather large ($>c$) probability
to be very small, of the order of $r^{2}$. However, since there is
also a probability $>c$ to escape $B(v,\frac{1}{2}N)$, this fact
has negligible impact on the moments of $t_{i+1}-t_{i}$.%
}, with its expectation being less or equal than\begin{equation}
U_{i}:=\begin{cases}
 Cr^{2} & i\textrm{ is odd}\\
 CN^{d}r^{2-d} & i\textrm{ is even}\end{cases}\label{eq:defUi}\end{equation}
 even after conditioning on the entry and exit points. In a formula,\begin{equation}
\mathbb{P}(t_{i+1}-t_{i}>\lambda U_{i}\, |\, R(u_{i})=y_{1}\textrm{ and }R(u_{i+1})=y_{2})\leq Ce^{-c\lambda }\label{eq:uiexpo}\end{equation}
for every $y_{1}\in B_{i}$ and $y_{2}\in B_{i+1}$ (see lemmas \ref{lem:condexpin}
and \ref{lem:condexpout}).

Define the variables \[
X_{i}:=(t_{i+1}-t_{i}-\mathbb{E}(t_{i+1}-t_{i}))/U_{0}\quad .\]
We wish to use lemma \ref{lem:pseudo} for the $X_{i}$'s. To get
(\ref{eq:Xiexpo}) we use (\ref{eq:uiexpo}) to see that $\mathbb{E}(t_{i+1}-t_{i})/U_{0}\leq CU_{i}/U_{0}\leq C$
and then use (\ref{eq:uiexpo}) again to get\begin{equation}
\mathbb{P}(X_{i}>\lambda \, |\, R(t_{i}),R(t_{i+1}))\leq Ce^{-c\lambda }\quad .\label{eq:XiRuiExpo}\end{equation}
Denote by $\mathcal{X}$ the event $X_{1},\dotsc ,X_{i-1},X_{i+1},\dotsc ,X_{n}$
and then\begin{align*}
\mathbb{P}(X_{i}>\lambda \, |\, \mathcal{X}) & =\mathbb{EP}(X_{i}>\lambda \, |\, R(t_{i}),R(t_{i+1}),\mathcal{X})=\mathbb{EP}(X_{i}>\lambda \, |\, R(t_{i}),R(t_{i+1}))\\
 & \leq \mathbb{E}Ce^{-c\lambda }=Ce^{-c\lambda }
\end{align*}
where the expectation above is with respect to $R(t_{i})$ and $R(t_{i+1})$.
This gives (\ref{eq:Xiexpo}).

The argument for (\ref{eq:Xipseudo}) requires the convergence of
the distributions. Start with the case of one $i$. Denote by $\mathcal{Y}$
the event $R(t_{i}),R(t_{i+1})$ and by $\mathcal{Z}$ the event $R(t_{i-\Delta }),\linebreak [0]R(t_{i+1+\Delta })$
for some $\Delta \in \{0,1,\dotsc \}$. Then\begin{align}
\mathbb{E}(X_{i}\, |\, \mathcal{Z}) & =\sum _{y\in B_{i}\times B_{i+1}}\mathbb{P}(\mathcal{Y}=y\, |\, \mathcal{Z})\cdot \mathbb{E}(X_{i}\, |\, \mathcal{Y}=y)\nonumber \\
 & =\sum _{y\in B_{i}\times B_{i+1}}\big (\mathbb{P}(\mathcal{Y}=y\, |\, \mathcal{Z})-\mathbb{P}(\mathcal{Y}=y)\big )\cdot \mathbb{E}(X_{i}\, |\, \mathcal{Y}=y)\label{eq:becauseEIS0}\\
 & \leq C\sum _{y\in B_{i}\times B_{i+1}}\big |\mathbb{P}(\mathcal{Y}=y\, |\, \mathcal{Z})-\mathbb{P}(\mathcal{Y}=y)\big |\label{eq:XiZ}
\end{align}
where the equality (\ref{eq:becauseEIS0}) is due to $\mathbb{E}X_{i}=0$.
Denote by $\pi _{k}(w,x)$ the probability to start from $w$ and
hit $x$ after $k$ moves of going from $B_{j}$ to $B_{j+1}$. In
a formula\[
\pi _{k}(w,x):=\mathbb{P}(R(u_{j+k})=x\, |\, R(u_{j})=w)\quad .\]
Of course, we mean that if $w\in \partial B(v,2r)$ then we take $j$
odd and in the opposite case we take $j>0$ even. Other than that
the value of $\pi _{k}$ is independent of $j$. With these notations
we get\begin{align*}
 & \mathbb{P}(\mathcal{Y}=(y_{1},y_{2}))=\mathbb{P}(R(t_{i})=y_{1})\pi _{1}(y_{1},y_{2})\\
 & \mathbb{P}(\mathcal{Y}=(y_{1},y_{2})\, |\, \mathcal{Z}=(z_{1},z_{2}))=\frac{\pi _{\Delta }(z_{1},y_{1})\pi _{1}(y_{1},y_{2})\pi _{\Delta }(y_{2},z_{2})}{\pi _{2\Delta +1}(z_{1},z_{2})}
\end{align*}
so\begin{align}
|\mathbb{P}(\mathcal{Y}=y\, |\, \mathcal{Z}=z)-\mathbb{P}(\mathcal{Y}=y)| & \leq \pi _{1}(y_{1},y_{2})\bigg (|\mathbb{P}(R(t_{i})=y_{1})-\pi _{\Delta }(z_{1},y_{1})|\; +\label{eq:P12}\\
 & +\; \left|\frac{\pi _{\Delta }(y_{2},z_{2})}{\pi _{2\Delta +1}(z_{1},z_{2})}-1\right|\pi _{\Delta }(z_{1},y_{1})\bigg )\quad .\nonumber 
\end{align}
Summing over $y$ the first half of (\ref{eq:P12}) we get \begin{eqnarray}
\lefteqn{\sum _{y_{1},y_{2}}\pi _{1}(y_{1},y_{2})|\mathbb{P}(R(t_{i})=y_{1})-\pi _{\Delta }(z_{1},y_{1})|=} &  & \label{eq:sum1}\\
 & \qquad \qquad  & =\sum _{y_{1}}|\mathbb{P}(R(t_{i})=y_{1})-\pi _{\Delta }(z_{1},y_{1})|\leq 2e^{-c\Delta }\nonumber 
\end{eqnarray}
where the last inequality is due to the exponential convergence of
the distributions in the form (\ref{eq:l1expo}) --- take $q_{1}$
to be the distribution of $R(t_{i-\Delta })$ and $q_{2}=\delta _{\{z_{1}\}}$
(the distance between any two distributions is always $\leq 2$).
For the second half of (\ref{eq:P12}), we use the form (\ref{eq:linfexpo})
for and get, under the assumption $\Delta >0$,\begin{align}
\lefteqn{\sum _{y_{1},y_{2}}\pi _{\Delta }(z_{1},y_{1})\pi _{1}(y_{1},y_{2})\left|\frac{\pi _{\Delta }(y_{2},z_{2})}{\pi _{2\Delta +1}(z_{1},z_{2})}-1\right|\leq } & \nonumber \\
 & \qquad \qquad \qquad \leq Ce^{-c\Delta }\sum _{y_{1},y_{2}}\pi _{\Delta }(z_{1},y_{1})\pi _{1}(y_{1},y_{2})=Ce^{-c\Delta }\quad .\label{eq:sum2}
\end{align}
 We used here (\ref{eq:linfexpo}) with $q_{1}=\delta _{\{y_{2}\}}$
and $q_{2}$ the distribution of $R(u_{j+\Delta +1})\, |\, R(u_{j})=z_{1}$
for the point $z_{2}$. Using (\ref{eq:sum1}), (\ref{eq:sum2}) and
(\ref{eq:P12}) in (\ref{eq:XiZ}) gives \begin{equation}
\mathbb{E}(X_{i}\, |\, R(t_{i-\Delta }),R(t_{i+1+\Delta }))\leq Ce^{-c\Delta }\quad .\label{eq:oneXi}\end{equation}
(the case $\Delta =0$ doesn't follow from the argumentation above,
but can be deduced, say, from (\ref{eq:XiRuiExpo})).

With (\ref{eq:oneXi}), proving (\ref{eq:Xipseudo}) is easy. Let
$i_{1},\dotsc ,i_{l}$ be some integers, all different, and let\[
\Delta _{j}=\bigg \lfloor \frac{1}{2}\bigg (\min _{\substack{ 1\leq m\leq l\\
 m\neq j}
}|i_{j}-i_{m}|-1\bigg )\bigg \rfloor \quad ,\]
so that the intervals $\left]i_{j}-\Delta _{j},i_{j}+1+\Delta _{j}\right[$
are disjoint. Let $\mathcal{X}$ be the event \[
R(t_{i_{1}-\Delta _{1}}),R(t_{i_{1}+1+\Delta _{1}}),\dotsc ,R(t_{i_{k}-\Delta _{k}}),R(t_{i_{k}+1+\Delta _{k}})\quad .\]
Then conditioning by $\mathcal{X}$ the events $X_{i}$ are independent
so we get\begin{align*}
\mathbb{E}(X_{i_{1}}\dotsm X_{i_{k}}\, |\, \mathcal{X}) & =\prod _{j=1}^{k}\mathbb{E}(X_{j}\, |\, \mathcal{X})=\prod _{j=1}^{k}\mathbb{E}(X_{j}\, |\, R(t_{i_{j}-\Delta _{j}}),R(t_{i_{j}+1+\Delta _{j}}))\\
 & \! \! \stackrel{(\ref {eq:oneXi})}{\leq }\prod _{j=1}^{k}Ce^{-c\Delta _{j}}\leq \prod _{j=1}^{k}C\exp (-c\min _{\substack{ 1\leq m\leq l\\
 m\neq j}
}|i_{j}-i_{m}|)
\end{align*}
 which immediately gives (\ref{eq:Xipseudo}) since\begin{align*}
\lefteqn{\mathbb{E}(X_{i_{1}}\dotsm X_{i_{k}}|\, X_{i_{k+1}},\dotsc ,X_{i_{l}})=} & \\
 & \qquad \qquad =\mathbb{E}\big (\mathbb{E}(X_{i_{1}}\dotsm X_{i_{k}}\, |\, \mathcal{X})\, \big |\, X_{i_{k+1}},\dotsc ,X_{i_{l}}\big )\leq \\
 & \qquad \qquad \leq \mathbb{E}\bigg (\prod _{j=1}^{k}C\exp (-c\min _{\substack{ 1\leq m\leq l\\
 m\neq j}
}|i_{j}-i_{m}|)\, \bigg |\, X_{i_{k+1}},\dotsc ,X_{i_{l}}\bigg )=\displaybreak [0]\\
 & \qquad \qquad =\prod _{j=1}^{k}C\exp (-c\min _{\substack{ 1\leq m\leq l\\
 m\neq j}
}|i_{j}-i_{m}|)
\end{align*}
with (\ref{eq:Xiexpo}) and (\ref{eq:Xipseudo}) established we can
invoke lemma \ref{lem:pseudo} and get \[
\mathbb{P}(|t_{n}-\mathbb{E}t_{n}|>\lambda U_{0})\leq Ce^{-c\lambda ^{2}}\quad .\]
Lemma \ref{lem:times} now follows since (\ref{eq:linfexpo}) shows
that $\mathbb{E}t_{2i+1}-t_{2i}$ converge exponentially to some $E_{\textrm{even}}$
and $\mathbb{E}t_{2i+2}-t_{2i+1}$ converge exponentially to some
$E_{\textrm{odd}}$ so \[
\left|\mathbb{E}t_{n}-n{\textstyle \frac{1}{2}}(E_{\textrm{even}}+E_{\textrm{odd}})\right|\leq CU_{0}\]
and for $\lambda <C\sqrt{n}$ this translation affects only the multiplicative
constant. Therefore taking $E=\frac{1}{2}(E_{\textrm{even}}+E_{\textrm{odd}})\approx N^{d}r^{2-d}$
and $\sigma =U_{0}\approx E$ we are done.
\end{proof}
\begin{lem}
\label{lem:abs}Let $b\in T_{N}^{d}$ and let $R$ be a random walk
on $T$ starting from $b$. Let $r<\frac{1}{8}N$, $v\in T_{N}^{d}$.
Let $t\in \mathbb{N}$ be some time. Then \begin{equation}
\mathbb{P}(\LE (R[0,t])\cap B(v,r)>\lambda f(r))\leq Ce^{-c\lambda }\quad \forall \lambda >0\label{eq:lemabs}\end{equation}
where $f$ is defined by (\ref{eq:deff}).
\end{lem}
\begin{proof}
Let $\lambda >0$ be some number. We note that we may assume $t<\lambda N^{d/2}$
since in time $\lambda N^{d/2}$ the probability to hit $b$ is $>1-Ce^{-c\lambda }$
and in this case the process starts afresh, memoryless. Let $t_{i}$
be stopping times defined by (\ref{eq:defti}). Let $E$ and $\sigma $
be defined by lemma \ref{lem:times} so that (\ref{eq:explamb2})
holds. 

The first case is $\lambda >C_{1}\log r$ for some $C_{1}$ sufficiently
large. This case is uninteresting for the following reason: lemma
\ref{lem:times} gives that for some $C_{2}$ sufficiently large,
if $n=\left\lfloor C_{2}\lambda r^{d-2}\right\rfloor $ then $\mathbb{P}(t_{n}\leq t)\leq Ce^{-c\lambda }$.
Let $k:=\max \{l:t_{l}\leq t\}$. If $k$ is even then $\LE (R[0,t])\cap B(v,r)\subset \LE (R[0,t_{k}])\cap B(v,r)$.
If $k$ is odd then \[
\#((\LE (R[0,t])\setminus \LE (R[0,t_{k}]))\cap B(v,r))\leq t_{k+1}-t_{k}\]
and this variable has the estimate (\ref{eq:tijump}) so it is uninteresting.
Therefore it is enough to calculate the loop-erased at the times $t_{k}$.
We get \begin{eqnarray*}
\lefteqn{\mathbb{P}(\#(\LE (R[0,t])\cap B(v,r))>\lambda f(r))\leq } &  & \\
 & \qquad  & \leq Ce^{-c\lambda }+\sum _{k=1}^{n}\mathbb{P}(\#(\LE (R[0,t_{k}])\cap B(v,r))>{\textstyle \frac{1}{2}}\lambda f(r))\leq \\
 &  & \leq Cne^{-c\lambda }\leq Cr^{d-2}\lambda e^{-c\lambda }\leq Cr^{d-2}e^{-c\lambda }\leq Cr^{d-2-cC_{1}}e^{-c\lambda }.
\end{eqnarray*}
(of course, all $c$'s in the last line are different). This shows
that for $C_{1}$ sufficiently large --- namely, $(d-2)/c$ where
$c$ is the last $c$ on the last line above, (\ref{eq:lemabs}) holds.
Thus this case is proved.

Therefore we shall assume that $\lambda <C\log r$. Let $n_{1}$ be
defined by \begin{align*}
n_{1}^{-} & :=\max \{n\textrm{ even}:t-nE>\lambda \sigma \sqrt{n}\}\\
n_{1}^{+} & :=\min \{n:t-nE<-\lambda \sigma \sqrt{n}\}
\end{align*}
Note that \[
n_{1}^{+}-n_{1}^{-}\leq C\lambda \sqrt{t/E}\leq CN^{d/2}\lambda ^{3/2}\sqrt{r^{d-2}}.\]
Let $\mathcal{E}_{1}$ be the event $|t_{n_{1}^{-}}-n_{1}^{-}E|>\lambda \sigma \sqrt{n_{1}^{-}}$
or $|t_{n_{1}^{+}}-n_{1}^{+}E|>\lambda \sigma \sqrt{n_{1}^{+}}$.
Lemma \ref{lem:times} gives us that $\mathbb{P}(\mathcal{E}_{1})<Ce^{-c\lambda }$.
We note that under $\neg \mathcal{E}_{1}$ we can {}``locate'' $t$,
$t_{n_{1}^{-}}<t<t_{n_{1}^{+}}$ and the interval is not very large,
$t_{n_{1}^{+}}-t_{n_{1}^{-}}<CN^{d}\lambda ^{3/2}\sqrt{r^{2-d}}$.
Let $\mathcal{E}_{2}$ be the event $\#\LE (R[0,t_{n_{1}}])\cap B(v,r))>\lambda f(r)$.
Lemma \ref{lem:mainpp} gives us that $\mathbb{P}(\mathcal{E}_{2})<Ce^{-c\lambda }$. 

We continue to define a short sequence of $n_{j}^{\pm }$ inductively:\begin{align*}
n_{j}^{-} & :=\max \{n\textrm{ even}:t-t_{N_{j-1}^{-}}-nE>\lambda \sigma \sqrt{n}\}\\
n_{j}^{+} & :=\min \{n:t-t_{N_{j-1}^{-}}-nE<-\lambda \sigma \sqrt{n}\}\\
N_{j}^{\pm } & :=n_{j}^{\pm }+\sum _{k=1}^{j-1}n_{k}^{-}
\end{align*}
Unlike $n_{1}^{\pm }$ which are just numbers, $n_{j}^{\pm }$, $j>1$
are events depending on $R[0,t_{N_{i-1}^{-}}]$. Define $\mathcal{E}_{2i-1}$
to be the event $|t_{N_{i}^{\pm }}-t|>\lambda \sigma \sqrt{n_{i}^{\pm }}$
(as before, we mean that either happens). Again, we get $\mathbb{P}(\mathcal{E}_{2i-1})<Ce^{-c\lambda }$.
Under $\neg (\mathcal{E}_{1}\cup \mathcal{E}_{3}\cup \dotsb \cup \mathcal{E}_{2i-3})$
we have\begin{equation}
n_{i}^{\pm }<C\frac{t-t_{N_{i-1}^{-}}}{E}\leq C\frac{t_{N_{i-1}^{+}}-t_{N_{i-1}^{-}}}{E}\stackrel{(\ref {eq:tnipm})}{\leq }C(i)\lambda ^{2-2^{-i+1}}r^{2^{-i+1}(d-2)}\label{eq:Nipm}\end{equation}
(the use of (\ref{eq:tnipm}) is inductively, for $i-1$). The addition
of $\neg \mathcal{E}_{2i-1}$ gives $t_{N_{i}^{-}}<t<t_{N_{i}^{+}}$
and\begin{equation}
t_{N_{i}^{+}}-t_{N_{i}^{-}}\leq \lambda \sigma \left(\sqrt{n_{i}^{+}}+\sqrt{n_{i}^{-}}\right)\stackrel{(\ref {eq:Nipm})}{\leq }C(i)N^{d}\lambda ^{2-2^{-i}}r^{(1-2^{-i})(2-d)}\quad .\label{eq:tnipm}\end{equation}
To use lemma \ref{lem:mainpp}, we need to define an auxiliary walk
$R'$,\[
R_{i}'(u)=R(u+t_{N_{i-1}^{-}})\quad R_{i}':\{-t_{N_{i-1}^{-}},\dotsc ,t-t_{N_{i-1}^{-}}\}\rightarrow T.\]
In other words, we consider the part of the walk until $t_{N_{i-1}^{-}}$
as fixed, and the part from $t_{N_{i-1}^{-}}$ to $t$ as the probabilistic
part. Of course, the stopping times $t_{j}'$ corresponding to $R_{i}'$
are simply $t_{j}'=t_{N_{i-1}^{-}+j}$. The fact that $N_{i}^{-}$
is even means that $R_{i}'(0)\not \in B(v,2r)$ and then lemma \ref{lem:mainpp}
will give that\begin{eqnarray}
\mathbb{P}(\mathcal{E}_{2i}) & \leq  & Ce^{-c\lambda }\nonumber \\
\mathcal{E}_{2i} & := & \{\#(L_{i}\cap B(v,r))>\lambda f(r)\}\label{eq:defE2i}\\
L_{i} & := & \LE ^{+}(R_{i}'[-t_{N_{i-1}^{-}},t_{N_{i}^{-}}-t_{N_{i-1}^{-}}])\label{eq:defLi}
\end{eqnarray}
 Note that we have now defined all the exceptional events $\mathcal{E}_{i}$:
the even ones are (\ref{eq:defE2i}) and the odd ones have been defined
slightly above.

When we said that the series $n_{i}^{\pm }$ is short, we meant that
we shall take it until $I$ defined by\[
I=\begin{cases}
 2 & d\geq 7\\
 3 & d=5,6\\
 C\log \epsilon ^{-1} & d=4\end{cases}\]
where $\epsilon $ is from (\ref{eq:deff}), which we consider as
a constant, so $I\leq C$. In particular\[
\mathbb{P}(\mathcal{E}_{1}\cup \dotsb \cup \mathcal{E}_{2i})\leq CIe^{-c\lambda }\leq Ce^{-c\lambda }.\]
The reason for this selection of $I$ is that with this $I$ it is
possible to do a simple estimate of the path between $t_{N_{I}^{-}}$
and $t$. For any $i$ we have \begin{align*}
N_{i}^{+}-N_{i}^{-} & <C\lambda \sqrt{(t-t_{N_{i-1}^{-}})/E}\leq C\lambda \sqrt{(t_{N_{i-1}^{+}}-t_{N_{i-1}^{-}})/E}\\
 & \! \! \stackrel{(\ref {eq:tnipm})}{\leq }C\lambda ^{2-2^{-i}}r^{2^{-i}(d-2)}\leq Cr^{2^{-i}(d-2)}\log ^{2}r
\end{align*}
(remember that $\lambda <C\log r$) and for $I$ this gives \[
N_{I}^{+}-N_{I}^{-}\leq Cf(r)r^{-2}\quad .\]
Therefore we may use (\ref{eq:tijump}) $N_{I}^{+}-N_{I}^{-}$ times,
to get \begin{equation}
\mathbb{P}\bigg (\sum _{\substack{ j=N_{I}^{-}\\
 j\textrm{ odd}}
}^{N_{I}^{+}}t_{j+1}-t_{j}>\lambda f(r)\bigg )\leq Ce^{-c\lambda }\label{eq:tNI_t}\end{equation}
which of course bounds also $\#(\LE (R[0,t])\cap B(v,r))-\#(\LE (R[0,t_{N_{I}^{-}}])\cap B(v,r))$.
Finally, the definitions of $R_{i}'$, $\LE $, $\LE ^{+}$ and $L_{i}$
(\ref{eq:defLi}) give \[
\LE (R[0,t_{N_{I}^{-}}])\subset L_{1}\cup L_{2}\cup \dotsb \cup L_{I}\]
and assuming $\neg (\mathcal{E}_{2}\cup \mathcal{E}_{4}\cup \dotsb \cup \mathcal{E}_{2i})$
we have from (\ref{eq:defE2i}) that \[
\#(\LE (R[0,t_{N_{I}^{-}}])\cap B(v,r))\leq I\lambda f(r)\leq C\lambda f(r)\]
and with (\ref{eq:tNI_t}) we finally get\[
\mathbb{P}(\#(\LE (R[0,t])\cap B(v,r))>\lambda f(r))\leq Ce^{-c\lambda }\quad .\]
and the lemma is proved.
\end{proof}
\begin{rem*}
By now the reader would not be surprised to learn that here too, if
one is willing to let go of a $\log $ factor then the proof gets
much simpler. Indeed, the arguments used for the case $\lambda >C\log r$
can be used for any $\lambda $ to get this result, and for this case
one does not need the precise estimates of lemma \ref{lem:times}
either, and the entire section may be reduced to half a page.
\end{rem*}
\begin{thm}
\label{thm:absolute}Let $b\in T_{N}^{d}$ and let $R$ be a random
walk on $T$ starting from $b$. Let $t\in \mathbb{N}$ be some time.
Then \[
\mathbb{P}(\LE (R[0,t])>\lambda f(N))\leq Ce^{-c\lambda }\quad \forall \lambda >0\]
where $f$ is defined by (\ref{eq:deff}).
\end{thm}
The theorem follows from lemma \ref{lem:abs} like theorem \ref{thm:upper}
follows from lemma \ref{lem:mainpp} (cover $T$ by balls etc.) and
we shall omit the proof.

\section{\label{sec:The-lower-bound}The lower bound}

We will use the concept of a cut time

\begin{defn*}
Let $R$ be a random walk on a graph, possibly with a stopping condition.
A time $t$ is called a \textbf{cut time} for $R$ if $R[0,t]\cap R\left]t,\infty \right[=\emptyset $.
\end{defn*}
Clearly, if $t$ is a cut time then $R(t)\in \LE (R)$. Further, all
$R(t_{i})$'s for different cut times $t_{i}$ are different. Therefore
it is possible to estimate the length of a loop-erased random walk
by counting cut times.

\begin{lem}
\label{lem:cuttime}Let $d\geq 5$. Let $R$ be a random walk on $T_{N}^{d}$
of length $L$ for some $L=\epsilon N^{d/2}$, $\epsilon $ sufficiently
small and $N>N_{0}(\epsilon )$. Let $X$ be the number of cut times
of $R$. Then\begin{equation}
\mathbb{E}X>cL\quad \mathbb{V}X<C\epsilon ^{2}L^{2}\quad .\label{eq:lemcuttime}\end{equation}

\end{lem}
As usual $\mathbb{V}$ denotes the variance, i.e.~$\mathbb{V}X:=\mathbb{E}X^{2}-(\mathbb{E}X)^{2}$.

\begin{proof}
Denote by $E_{t}$ the event that $t$ is a cut time. Easily,\[
1-\mathbb{P}(E_{t})\leq \sum _{s_{1}=0}^{t}\sum _{s_{2}=t+1}^{L}\mathbb{P}(R(s_{1})=R(s_{2}))\quad .\]
Now for $|s_{i}-t|\leq N$ this is identical to the equivalent problem
on $\mathbb{Z}^{d}$ which is well known (see \cite{L96}) so we get
\[
\mathbb{P}(R[\max 0,t-N,t]\cap R\left]t,\min L,t+N\right]\neq \emptyset )<1-c\quad .\]
For other $s_{i}$ we use the easy \begin{equation}
\mathbb{P}(R(s_{1})=R(s_{2}))\leq C\min \{N^{2},|s_{1}-s_{2}|\}^{-d/2}\label{eq:Rs1Rs2simple}\end{equation}
 to get \[
1-\mathbb{P}(E_{t})<1-c+C\epsilon ^{2}+CN^{2-d/2}\]
therefore for $\epsilon $ sufficiently small and $N$ sufficiently
large we get $\mathbb{P}(E_{t})>c$ which gives the first part of
(\ref{eq:lemcuttime}) --- $\mathbb{E}X>cL$. For the second part,
we examine the covariance of $E_{t_{1}}$ and $E_{t_{2}}$ for some
$t_{1}<t_{2}$. Denote $t=\left\lfloor \frac{1}{2}(t_{1}+t_{2})\right\rfloor $
and \[
E_{1}'=\mathbb{P}(R[0,t_{1}]\cap \left]t_{1},t\right]=\emptyset )\quad E_{2}'=\mathbb{P}(R[t,t_{2}]\cap R\left]t_{2},L\right]=\emptyset )\quad .\]
We note that $E_{1}'$ and $E_{2}'$ are independent and therefore
$\cov E_{1}',E_{2}'=0$. On the other hand, summing (\ref{eq:Rs1Rs2simple})
we get\begin{align*}
|\mathbb{P}(E_{1}')-\mathbb{P}(E_{t_{1}})| & \leq \mathbb{P}(R[0,t_{1}]\cap R[t,L]\neq \emptyset )\leq \sum _{s_{1}=0}^{t_{1}}\sum _{s_{2}=t}^{L}C\min \{N^{2},|s_{2}-s_{1}|\}^{-d/2}\\
 & \leq C\sum _{s_{1}=0}^{t_{1}}|t-s|^{1-d/2}+\epsilon N^{-d/2}\leq C(|t_{2}-t_{1}|^{2-d/2}+\epsilon ^{2})
\end{align*}
so we get the same for the covariance of $E_{t_{i}}$,\[
\cov E_{t_{1}},E_{t_{2}}\leq C|t_{2}-t_{1}|^{2-d/2}+C\epsilon ^{2}\quad .\]
Summing these for all $t_{i}$'s we get the second half of the lemma.
\end{proof}
\begin{lem}
\label{lem:LEtbig}Let $d\geq 5$. Let $b\in T_{N}^{d}$ and let $R$
be a random walk on $T$ starting from $b$. Let $t\in \mathbb{N}$,
$t>N^{d/2}$ and $\lambda >N^{-1/2}$. Then \[
\mathbb{P}(\#\LE (R[0,t])\leq \lambda N^{d/2})\leq C\lambda \]

\end{lem}
\begin{proof}
We may assume without loss of generality that $\lambda \leq c$ for
some constant. Let $C_{1}$ be some constant which will be fixed later.
Define\[
u:=t-C_{1}\lambda N^{d/2}.\]
(we assume here $\lambda <1/C_{1}$, as we may). Denote by $X$ the
number of cut times in the segment $[u,t]$. Lemma \ref{lem:cuttime}
shows that $\mathbb{E}X>c(t-u)=cC_{1}\lambda N^{d/2}$. Pick $C_{1}$
sufficiently large such that $\mathbb{E}X>3\lambda N^{d/2}$. Lemma
\ref{lem:cuttime} also gives $\mathbb{V}X\leq C\lambda ^{4}N^{d/2}$
and then\[
\mathbb{P}(X\leq 2\lambda N^{d/2})\leq C\lambda ^{2}.\]
Next we want to estimate\[
\mathbb{P}(\LE (R[0,u])\cap R[u+N^{2},t]\neq \emptyset .)\]
Define $Y=\#\{\LE (R[0,u])\cap R[u+N^{2},t]\}$. If we assume $\#\LE (R[0,u])\leq \mu N^{d/2}$,
then because $R(u+N^{2}$) is distributed $\approx $ uniformly on
$T$ we get \[
\mathbb{E}(Y\, |\, \#\LE (R[0,u])\leq \mu N^{d/2})\approx N^{-d}(\#\LE (R[0,u]))(t-u-N^{2})\approx \mu \lambda \]
(this is the only place we use the assumption $\lambda >N^{-1/2}$).
Without the assumption $\#\LE (R[0,u])\leq \mu N^{d/2}$ we get\begin{align*}
\mathbb{E}Y & \leq \sum _{\mu =0}^{\infty }\mathbb{P}(\#\LE (R[0,u])\leq \mu N^{d/2})\cdot \mathbb{E}(Y\, |\, \#\LE (R[0,u])\leq (\mu +1)N^{d/2})\\
 & \leq \sum _{\mu =0}^{\infty }Ce^{-c\mu }\mu \lambda \leq C\lambda 
\end{align*}
 and hence $\mathbb{P}(Y>0)\leq C\lambda $. Under the assumption
$Y=0$ every cut point of $R[u,t]$ above $u+N^{2}$ is in $\LE (R[0,t])$
and the lemma follows.
\end{proof}

\begin{proof}
[Proof of theorem \ref{thm:dge5}]Let $R'$ be a random walk starting
from $b$ with no stopping condition. Define events\begin{align*}
\mathcal{X}(v,t) & =\{R'(t)=v\wedge v\not \in R'\left[0,t\right[\}\\
\mathcal{Y}(v,t) & =\{R'(t)=v\wedge \#\LE (R[0,t])\leq \lambda N^{d/2}\}.
\end{align*}
Now, $\sum _{v}\mathbb{P}(\mathcal{X}(v,t))$ is simply the probability
that a random walk reaches its end point for the first time, or equivalently
by symmetry, the probability that it never returned to its starting
point, therefore it is easy to calculate\[
\sum _{v\in T}\mathbb{P}(\mathcal{X}(v,t))\leq Ce^{-ctN^{-d}}\quad \forall t.\]
Next, for $t>N^{d/2}$, lemma \ref{lem:LEtbig} gives\[
\sum _{v\in T}\mathbb{P}(\mathcal{Y}(v,t))\leq C\lambda \quad \forall t>N^{d/2}.\]
Finally, note that \[
\sum _{t=0}^{\infty }\mathbb{P}(\mathcal{X}(v,t))=1\quad \forall v\in T.\]
With these three facts we get, for any parameter $\mu >0$,\begin{eqnarray*}
\lefteqn{\sum _{t,v}\mathbb{P}(\mathcal{X}(v,t)\setminus \mathcal{Y}(v,t))\geq } &  & \\
 & \qquad  & \geq N^{d}-\left(\sum _{t=0}^{N^{d/2}}+\sum _{t=\mu N^{d}}^{\infty }\right)\sum _{v\in T}\mathcal{P}(\mathcal{X}(v,t))-\sum _{t=N^{d/2}}^{\mu N^{d}}\sum _{v\in T}\mathbb{P}(\mathcal{Y}(v,t))\geq \\
 &  & \geq N^{d}(1-Ce^{-c\mu }-N^{-d/2}-C\mu \lambda ).
\end{eqnarray*}
 Picking $\mu =C\log \lambda $ for some $C$ sufficiently large will
prove the theorem.
\end{proof}

\subsection{Remarks on alternative approaches}

The first alternative approach to the proof of the lower bound is
as follows: prove a conditioned version of lemma \ref{lem:cuttime},
namely

\begin{lem*}
Let $b$ and $e$ be two points on $T_{N}^{d}$ with $|b-e|>cN$.
Let $R$ be a random walk on $T_{N}^{d}$ of length $L$ for some
$L=\epsilon N^{d/2}$, $\epsilon $ sufficiently small starting from
$b$ and conditioned to end at $e$. Let $X$ be the number of cut
times of $R$. Then\[
\mathbb{E}X>cL\quad \mathbb{V}X<C\epsilon ^{2}L^{2}\quad .\]

\end{lem*}
This lemma allows to prove a version of theorem \ref{thm:dge5} for
any points far enough, not just two random points. Further, it allows
to avoid the need to use absolute times, and just work directly with
the times $t_{i}$ for some arbitrary ball. In other words, to show
that the loop-erased random walk from $b$ to $e$ is long with high
probability, define an arbitrary ball $B$, show that at the stopping
times $t_{i}$ corresponding to $B$ the entire loop-erased random
walk is quite small (this is quite simple) and then show that the
random walk from the last $t_{i}$ to $e$ has many cut points using
the lemma above. 

The proof of this lemma requires no new ideas when compared with lemma
\ref{lem:cuttime}. However, it is very technical, and quite long,
which is the main reason we chose the approach above. In some sense
we do not consider the length of section \ref{sec:Absolute-times}
as an indication that the approach we chose is more complicated because
the result (theorem \ref{thm:absolute}) is trivial if one can afford
to lose a $\log $ factor (and also because the result is quite natural).

Another approach is the use of the uniform spanning tree and Wilson's
algorithm (see \cite{W96}). Roughly, one might hope to show that
the loop-erased random walk is long by constructing an appropriate
partial UST, and then showing that the random walk $R$ starting from
some point $b$ and stopped on the partial UST is not too long (therefore
no complicated self interactions, as in lemma \ref{lem:cuttime})
and not too short, so $\LE (R)$ can be proved to be long. Since the
loop-erased random walk from $b$ to some other point $e$ (say inside
the partial UST) contains $\LE (R)$, this will be enough. Alternatively,
one can take two random walks $R$ and $R'$ starting from $b$ and
$e$ respectively and stopped on the partially constructed UST, and
calculate the probabilities that at least one is long and that they
do not intersect. Both approaches allow to generalize theorem \ref{thm:dge5}
from a random end point to any end point (naturally, if $b$ and $e$
are very close then with positive probability the loop-erased random
walk from $b$ to $e$ is short. However, one can show that there
is a positive probability for the loop-erased random walk to be long,
i.e.~$\approx N^{d/2}$). 

A third strategy using the UST is as follows. Notice that the harmonic
measure on a partially constructed UST is roughly uniform --- this
follows since the escape probabilities from a typical small ball are
positive. If one wants to estimate the probability that the loop-erased
random walk between $b$ and $e$ is $\leq \lambda N^{d/2}$, construct
a partial UST containing $b$ up until its size is $\approx (1/\lambda )N^{d/2}$,
and then estimate that the number of vertices in the tree with distant
$\leq \lambda N^{d/2}$ from $b$ is $\approx \lambda N^{d/2}$ ,
so the harmonic measure is $\approx \lambda ^{2}$. This approach
gives (in addition to the fact that $e$ may be arbitrary) stronger
estimates than the $\lambda \log \lambda ^{-1}$ of theorem \ref{thm:dge5}
--- formalizing these arguments we were able to show $\mathbb{P}(\#\LE (R)\leq \lambda N^{d/2})\leq C\lambda ^{2}\log \lambda ^{-1}$,
and we believe that the true value is, as in the case of the complete
graph, $\lambda ^{2}$. The difference between $\lambda \log \lambda ^{-1}$
and $\lambda ^{2}\log \lambda ^{-1}$ is significant in the following
sense: the weaker estimate does not prove that the UST has a true
branching nature: even points that are distributed linearly along
a path of length $N^{d/2}$ satisfy the requirement that $\mathbb{P}(\LE \leq \lambda N^{d/2})\leq C\lambda $.
However, the estimate $\lambda ^{2}\log \lambda ^{-1}$ allows to
deduce non-trivial facts about the branching structure of the UST.

None of these methods work in dimension $4$, and the culprit is always
the same: in dimension $4$ our methods do not show that within the
mixing time the probability of hitting the loop-erased random walk
is small. In other words, to get a lower bound for dimension $4$
one must either show a very precise upper estimate (not much different
from the conjectured precise value) or alternatively show indirectly
that the mixing time is smaller than the hitting time of the loop-erased
random walk.

\appendix

\section{Proofs of known and unsurprising facts}

\renewcommand{\thelem}{A.\arabic{lem}}
\setcounter{lem}{0}The harmonic potential on $\mathbb{Z}^{d}$, $d>2$, is the unique
bounded function $a$ satisfying \[
\Delta a(z)=\begin{cases}
 1 & z=\vec{0}\\
 0 & \textrm{otherwise}\end{cases}\]
 and $a(\infty )=0$ where $\Delta $ stands for the discrete Laplacian.
It is well known (see e.g.~\cite[theorem 1.5.4]{L96} %
\footnote{\cite{L96} only shows $a(v)=\alpha |v|^{2-d}+O(|v|^{\epsilon -d})$,
but this is completely sufficient for our purposes.%
} or \cite[theorem 5]{KS}) that $a(v)=\alpha |v|^{2-d}+O(|v|^{-d})$.

\begin{lem}
\label{lem:BinBpc}Let $B_{1}=B(x_{1},r_{1})\subset B_{2}=B(x_{2},r_{2})\subset T_{N}^{d}$,
$r_{2}\leq C_{1}r_{1}$. Let $v\in B_{2}\setminus B_{1}$ satisfy
$d(v,\partial B_{2})\geq c_{1}r_{1}$. Let $R$ be a random walk starting
from $v$ and stopped on $\partial B_{1}\cup \partial B_{2}$. Let
$p$ be the probability that $R$ hits $\partial B_{1}$. Then $p\geq c(c_{1},C_{1})$.
\end{lem}
We assume here that a ball (e.g.~$B_{2}$) satisfies $r_{2}<\frac{1}{2}N$
i.e.~it does not wrap itself because we are on a torus. This assumption
holds for all balls in this appendix, and we will not repeat it.

\begin{proof}
Clearly, we may assume $r_{1}$ is sufficiently large in the sense
that $r_{1}>C(c_{1},C_{1})$. Since we are dealing with a process
completely inside $B_{2}$, we may assume we are in $\mathbb{Z}^{d}$.
Assume first that $|x_{1}-v|<\frac{1}{2}d(x_{1},\partial B_{2})$.
Since $a_{1}(v):=a(v-x_{1})$ is harmonic on $B_{2}\setminus B_{1}$,
$a_{1}(R)$ is a martingale, and if we define $\tau $ to be the stopping
time on $\partial B_{1}\cup \partial B_{2}$ then we get $a_{1}(v)=\mathbb{E}a_{1}(R(\tau ))$,
so\begin{align}
a_{1}(v) & =p\mathbb{E}(a_{1}(R(\tau ))\, |\, R(\tau )\in \partial B_{1})+(1-p)\mathbb{E}(a_{1}(R(\tau ))\, |\, R(\tau )\in \partial B_{2})\nonumber \\
 & \leq p\alpha r_{1}^{2-d}(1+o(1))+(1-p)\alpha d(x_{1},\partial B_{2})^{2-d}(1+o(1))\label{eq:hp1}
\end{align}
(the $o(1)$ notations are as $r_{1}\rightarrow \infty $ and may
depend on $c_{1}$ and $C_{1}$) and from $a_{1}(v)=\alpha |x_{1}-v|^{2-d}(1+o(1))$
we get\begin{equation}
p\geq \frac{|x_{1}-v|^{2-d}-d(x_{1},\partial B_{2})^{2-d}}{r_{1}^{2-d}-d(x_{1},\partial B_{2})^{2-d}}(1+o(1))\geq c\label{eq:hp2}\end{equation}
for $r_{1}$ sufficiently large. In the case $|x_{1}-v|\geq \frac{1}{2}d(x_{1},\partial B_{2})$,
we can find a sequence of balls $B(y_{n},s_{n})$ of length $\leq C(c_{1},C_{1})$
and each $s_{n}\geq c(c_{1},C_{1})r_{1}$ such that $|y_{1}-v|\leq \frac{1}{2}d(y_{1},\partial B_{2})$
and $\forall w\in \partial B(y_{i},s_{i})$, $|y_{i+1}-w|\leq \frac{1}{2}d(y_{i+1},\partial B_{2})$.
Notice that this is possible because $d(v,\partial B_{2})\geq c_{1}r_{1}$.
The previous case now gives that the probability that the random walk,
after hitting $B(y_{i},s_{i})$ will continue to $B(y_{i+1},s_{i+1})$
is $\geq c_{2}(c_{1},C_{1})$. Since it needs to perform only $C(c_{1},C_{1})$
such steps in order to hit $B_{1}$, we get $p\geq c_{2}^{C}=c_{3}(c_{1},C_{1})$.
\end{proof}
\begin{lem}
\label{lem:BoutBpc}Let $B(x_{1},r_{1}),B(x_{2},r_{2})$ be two balls,
$r_{2}\leq C_{1}r_{1}$ and $|x_{1}-x_{2}|\leq C_{1}r_{1}$; and let
$v\not \in B_{2}\cup B_{1}$ satisfy $d(v,B_{2})\geq c_{1}r_{2}$.
Then $p\geq c(c_{1},C_{1})$ where $p$ is as above.
\end{lem}
\begin{proof}
Assume first that $d(v,B_{1})\leq \frac{1}{2}d(B_{2},B_{1})$ where
$d(B_{1},B_{2})$ stands for the distance between the two balls in
the usual sense. Let $B_{i}'$ be the sets $B_{i}$ considered as
subsets of $\mathbb{Z}^{d}$ and let $S_{i}=B_{i}'+N\mathbb{Z}^{d}$
i.e.~$S_{i}$ is the preimage of $B_{i}$ by the quotient map $\mathbb{Z}^{d}\to T_{N}^{d}$.
Let $R'$ be a simple random walk on $\mathbb{Z}^{d}$ starting from
$v$ (we consider $v$ and the $B_{i}$'s as subsets of $\mathbb{Z}^{d}$
as well, say by locating them in $[0,N]^{d}$). Then\begin{align*}
p & =\mathbb{P}(R'\textrm{ hits }S_{1}\textrm{ before }S_{2})\geq \mathbb{P}(R'\textrm{ hits }B_{1}\textrm{ before }S_{2})\geq \\
 & \geq \mathbb{P}(R'\textrm{ hits }B_{1}\textrm{ before }\partial B(x_{1},r_{1}+d(B_{1},B_{2})))\stackrel{(*)}{\geq }c
\end{align*}
where $(*)$ comes from the same harmonic potential arguments as (\ref{eq:hp1})-(\ref{eq:hp2}).

If $d(v,B_{1})>\frac{1}{2}d(B_{2},B_{1})$ but we have both \begin{align*}
d(v,B_{1}) & \leq (2C_{1}+2)r_{1}\\
d(v,B_{2}) & \geq c_{1}r_{1}
\end{align*}
then the same ball-sequence argument as in the previous lemma gives
$p\geq c$;

If $d(v,B_{1})>(2C_{1}+2)r_{1}$, let $\tau $ be the hitting time
of $B_{3}:=B(x_{1},(2C_{1}+2)r_{1})$, then \begin{equation}
p=\mathbb{EP}(R'\textrm{ starting from }R(\tau )\textrm{ hits }B_{1}\textrm{ before }B_{2})\geq \mathbb{E}c=c\label{eq:contB3}\end{equation}
where here $R'$ is a simple random walk on $T_{N}^{d}$ (differing
from $R$ only by the starting point), the expectation $\mathbb{E}$
is over the distribution of $R(\tau )$ and the inequality comes from
the previous two cases.

Finally, if $d(v,B_{2})<c_{1}r_{1}$ define $\tau $ the hitting time
of $B_{4}:=B(x_{2},c_{1}r_{1})$. The harmonic potential at $x_{2}$
with a calculation similar to (\ref{eq:hp1})-(\ref{eq:hp2}) shows
that the probability to hit $B_{4}$ before $B_{2}$ is $\geq c$.
After hitting $B_{4}$ a calculation similar to (\ref{eq:contB3})
gives that $p\geq c$.
\end{proof}
\begin{lem}
\label{lem:poutB}Let $d(B(x_{1},r_{1}),B(x_{2},r_{2}))\geq c_{1}r_{1}$,
$r_{2}\geq c_{1}r_{1}$ and $|x_{1}-x_{2}|\leq C_{1}r_{2}$; and let
$v\in \partial B_{1}$. Let $R$ be a random walk starting from $v$
and stopped on $\partial B_{2}\cup \{x_{1}\}$. Let $p$ be the probability
that $R$ hits $x_{1}$. Then $p\approx r_{1}^{2-d}$.
\end{lem}
In the formulation of the lemma, and in its proof, all constants implicit
in the $\approx $ signs might depend on $c_{1}$ and $C_{1}$.

\begin{proof}
Let $B_{3}=B(x_{1},\frac{1}{2}r_{1})$. Define stopping times $t_{i}$
similarly to (\ref{eq:defti}), as follows: $t_{0}:=0$ and \begin{align*}
t_{2i+1} & :=\{t>t_{2i}:R(t)\in \partial B_{3})\\
t_{2i} & :=\{t>t_{2i-1}:R(t)\in \partial B_{1}\cup \{x_{1}\}\}.
\end{align*}
 Define also $\tau $ the hitting time of $\partial B_{2}$. The usual
harmonic potential calculations (use the harmonic potential around
$x_{1}$) show that the probability of a random walk starting from
any $v\in \partial B_{3}$ to hit $x_{1}$ before exiting $B_{1}$
is $\approx r_{1}^{2-d}$. Hence, $\mathbb{P}(R(t_{2i})=x_{1}|R(t_{2i-1}))\leq Cr_{1}^{2-d}$.
Lemma \ref{lem:BoutBpc} shows that a random walk starting from any
point in $\partial B_{1}$ has probability $\geq c$ to hit $B_{2}$
before hitting $B_{3}$. Therefore, the probability to get $t_{i}>\tau $
decreases exponentially in $i$ i.e.~$\mathbb{P}(t_{2i}<\tau )\leq Ce^{-ci}$
and hence\[
\mathbb{P}(R(t_{2i})=x_{1}\textrm{ and }t_{2i}<\tau )\leq Ce^{-ci}r_{1}^{2-d}\]
so\[
p\leq \sum _{i=0}^{\infty }\mathbb{P}(R(t_{2i+1})=x_{1}\textrm{ and }t_{2i+1}<\tau )\leq Cr_{1}^{2-d}.\]
The inequality $p\geq cr_{1}^{2-d}$ follows easily from (\ref{eq:hp1})-(\ref{eq:hp2})
for the harmonic potential at $x_{1}$ and the requirement $d(B_{1},B_{2})\geq c_{1}r_{1}$.
\end{proof}
\begin{lem}
\label{lem:hitp}Let $|v|<(1-c_{1})r$ and $w\in \partial B(0,r)$.
Then the probability $p$ that a random walk starting from $v$ will
exit $B$ in $w$ is $\approx r^{1-d}$. Without the restriction $|v|<(1-c_{1})r$
one has $p\leq C(r-|v|)^{1-d}$
\end{lem}
In the formulation of the lemma, and in its proof, all constants implicit
in the $\approx $ signs might depend on $c_{1}$.

\begin{proof}
In $\mathbb{Z}^{d}$, $d>2$, the probability of a walk starting from
$v$ to never return is $>c$. Hence the probability to hit $\partial B(0,r)$
before returning to $v$ is $>c$, and this event is identical on
$\mathbb{Z}^{d}$ and on the torus. The symmetry of the random walk
shows that $p$ is $\approx $ to the probability that a random walk
starting from $w$ will hit $\partial B\cup \{v\}$ in $v$ (the quotient
is exactly the probability of a random walk starting from $v$ to
return to $v$ before exiting $\partial B$, which, as we just discussed,
is $\approx 1$)%
\footnote{When we say {}``hit'' we mean at time $>0$, so that these probabilities
are not simply $1$.%
}. This probability can be calculated in three steps as follows. First,
the probability of a random walk starting from $w$ to hit $\partial B(0,\frac{2}{3}r+\frac{1}{3}|v|)$
is $\approx (r-|v|)^{-1}$: this uses an argument similar (\ref{eq:hp1})-(\ref{eq:hp2})
using the harmonic potential $a$ at $0$, but here we need the precise
estimate $a(x)=|x|^{2-d}+O(|x|^{-d})$ or at least $a(x)=|x|^{2-d}+O(|x|^{1-d})$
(see, e.g.~\cite[lemma 3]{K87} for a detailed version of this calculation).
Next, if $|v|<r(1-c_{1})$, use lemma \ref{lem:BinBpc} to show that
continuing from any point on $\partial B(0,\frac{2}{3}r+\frac{1}{3}|v|)$
the probability to hit $B(v,\frac{1}{3}(r-|v|))$ is $\approx 1$
--- if $|v|\geq r(1-c_{1})$, we only estimate that this probability
is $\leq 1$. Finally, the same (\ref{eq:hp1})-(\ref{eq:hp2}) argument
with the harmonic potential at $v$ shows that starting from any point
on $\partial B(v,\frac{1}{3}(r-|v|))$, the probability to hit $v$
before hitting $\partial B$ is $\approx (r-|v|)^{2-d}$.
\end{proof}
A similar calculation works when $(1+c_{1})r\leq ||v||$ and $p$
is the probability the random walk will hit $B$ in $w$, using lemma
\ref{lem:BoutBpc} instead of lemma \ref{lem:BinBpc} and lemma \ref{lem:poutB}
in the third step.

\begin{lem}
\label{lem:cond}Let $v\in B(0,r)\subset B(0,2r)\subset T_{N}^{d}$.
Let $R$ be a random walk starting from $v$ and stopped on $\partial B(0,2r)$.
Let $R_{x}$ be a random walk starting from $v$ and conditioned to
hit $\partial B(0,2r)$ at a specific point $x$. Then $R\cap B(0,r)\approx R_{x}\cap B(0,r)$
where $\approx $ means that the probabilities of any event are equal
up to a constant.
\end{lem}
\begin{proof}
Let $t$ be the last time when $R(t)\in B(0,r)$. Let $w=R(t)$. For
any $w$, the probability of an (unconditioned) random walk starting
from $w$ to hit $B(0,r)\cup \partial B(0,2r)$ in $\partial B(0,2r)$
is $\approx r^{-1}$. The probability to hit $x$ is $\approx r^{-d}$.
This independence from $w$ finishes the lemma. Both estimates are
easily proved as in the previous lemma.
\end{proof}
\begin{lem}
\label{lem:rrw_smooth}Let $R$ be a random walk starting from $v\in B(0,r)$
and let $w\in B(0,r)$. Let $t>r^{2}$. Let $p$ be the probability
that $R[0,t]\subset B(0,r)$ and $R(t)=w$. Then\[
p\leq Cr^{-d}e^{-ctr^{-2}}\quad .\]

\end{lem}
\begin{proof}
Starting from any $v\in B(0,r)$, after $r^{2}$ steps the random
walk has probability $>c$ to exit $B(0,r)$. This shows, clearly,
that the probability that $R[0,t-r^{2}]\subset B(0,r)$ is $\leq Ce^{-ctr^{-2}}$.
For any $x\in B(0,r)$, the probability that a random walk $R'$ starting
from $x$ satisfies $R'(r^{2})=w$ is $\leq Cr^{-d}$.
\end{proof}
\begin{lem}
\label{lem:erw}Let $v\in B(0,r)$ and let $R$ be a random walk starting
from $v$ and stopped on $\partial B(0,r)$. Let $w\in \partial B(0,r)$.
Then the probability $p$ that $R$ hits $w$ satisfies\[
p\leq C_{1}|v-w|^{1-d}\quad .\]

\end{lem}
\begin{proof}
Denote $s=|v-w|$. We shall prove the lemma using an induction process
that assumes the lemma holds for $1,\dotsc ,\frac{1}{2}s$ and proves
it for $\frac{1}{2}s+1,\dotsc ,s$.

Denote $d(v,\partial B(0,r))=\epsilon s$. The first thing to note
is that the lemma holds if $\epsilon >c$ with no need for induction
(in the sense that $p\leq C_{2}(c)|v-w|^{1-d}$), due to the second
part of lemma \ref{lem:hitp}.

It is for the case of small $\epsilon $ that we need the induction
process. Let $\delta =c_{1}\log ^{-1}\epsilon ^{-1}$ for some $c_{1}$
which will be fixed later. Let \begin{align*}
D & :=B(0,r)\setminus B(0,r-\delta s)\\
E: & =D\cap B(v,{\textstyle \frac{1}{2}}s)\quad .
\end{align*}
Examine the exit probabilities of $R$ from $E$. Let $p_{2}$ be
the probability that $R$ exits $E$ at $\partial B(0,r-\delta s)\cap \partial E$.
Then\[
p_{2}\leq \mathbb{P}(R\textrm{ exits }D\textrm{ at }\partial B(0,r-\delta s))\leq C\epsilon /\delta \]
where the second inequality comes from the harmonic potential at zero.
Let $p_{3}$ be the probability that $R$ exits $E$ at $\partial B(v,\frac{1}{2}s)\cap \partial E$.
It is easy to see that for any $x\in D$, the probability that $R$
exits $B(x,2\delta s)$ without hitting $\partial D$ is $<1-c$.
Therefore \[
p_{3}<(1-c)^{(s/2)/(2\delta s)}=e^{-cc_{1}^{-1}\log \epsilon ^{-1}}\quad .\]
Therefore for $c_{1}$ sufficiently small we would get $p_{3}\leq \epsilon $.
Together these two give \[
p_{2}+p_{3}\leq C\epsilon \log \epsilon ^{-1}\quad .\]
Since $E\cap B(w,\frac{1}{2}s)=\emptyset $ we get that the probability
of $R$ to hit $\partial B(w,\frac{1}{2}s)$ before exiting $B(0,r)$
is $\leq C\epsilon \log \epsilon ^{-1}$. The induction assumption
gives that the probability to hit $w$ after hitting $\partial B(w,\frac{1}{2}s)$
is $\leq C_{1}(\frac{1}{2}s)^{1-d}$. Therefore\[
p\leq C_{1}s^{1-d}\cdot (C2^{d-1}\epsilon \log \epsilon ^{-1})\quad .\]
For $\epsilon <c_{2}$ this will be $\leq C_{1}s^{1-d}$ and this
case is finished too. The lemma is now finished because for $\epsilon <c_{2}$
the induction process works for any $C_{1}$, and for $\epsilon \geq c_{2}$
the first case allows to define $C_{1}:=C_{2}(c_{2})$.
\end{proof}
A similar calculation works when $v$ is outside $B(0,r)$, and the
random walk is stopped when hitting $B(0,r)$ and the conclusion is
$p\leq C\min \{|v-w|,r\}^{1-d}$. 

\begin{lem}
\label{lem:condexpin}Let $v\in B(0,\frac{1}{2}r)$ and $w\in \partial B(0,r)$.
Let $R$ be a random walk started from $v$ and conditioned to exit
$B(0,r)$ in $w$. Let $t$ be the exit time. Then\[
\mathbb{P}(t>\lambda r^{2})\leq Ce^{-c\lambda }\quad .\]

\end{lem}
\begin{proof}
We may assume $\lambda >1$. Let $R'$ be an unconditioned walk starting
from $v$. Let \[
A_{n}=(B(w,2^{n})\setminus B(w,2^{n-1}))\cap B(0,r)\quad .\]
Lemma \ref{lem:rrw_smooth} shows that \[
\mathbb{P}(R'[0,\lambda r^{2}]\subset B(0,r)\textrm{ and }R'(\lambda r^{2})\in A_{n})\leq C\#A_{n}r^{-d}e^{-c\lambda }\quad .\]
Lemma \ref{lem:erw} shows that for every $x\in A_{n}$, the probability
of a random walk starting at $x$ to exit $B(0,r)$ at $w$ is $\leq C2^{n(1-d)}$.
Therefore by lemma \ref{lem:condexpin}\begin{align*}
\lefteqn{\mathbb{P}(R'[0,\lambda r^{2}]\subset B(0,r)\textrm{ and }R'(\lambda r^{2})\in A_{n}\textrm{ and }R'\textrm{ hits }w)\leq } & \\
 & \qquad \qquad \qquad \qquad \leq C(\#A_{n})2^{n(1-d)}r^{-d}e^{-c\lambda }\leq Cr^{-d}2^{n}e^{-c\lambda }
\end{align*}
and we get\[
\mathbb{P}(R'[0,\lambda r^{2}]\subset B(0,r)\textrm{ and }R'\textrm{ hits }w)\leq Cr^{1-d}e^{-c\lambda }\sum _{n=1}^{\left\lfloor \log r\right\rfloor }r^{-1}2^{n}\leq Cr^{1-d}e^{-c\lambda }\quad .\]
Since the probability of $R'$ to hit $w$ is $>cr^{1-d}$ (by lemma
\ref{lem:hitp}), we are done.
\end{proof}
\begin{lem}
\label{lem:sumexpexp}Let $X_{i}$ be events with a past-independent
exponential estimate, namely\[
\mathbb{P}(X_{i}>\lambda E\, |\, X_{i-1},\dotsc ,X_{1})\leq C_{1}e^{-c_{1}\lambda }.\]
Then\[
\mathbb{P}\left(\sum _{i=1}^{n}X_{i}>\lambda nE\right)\leq Ce^{-c\lambda }.\]

\end{lem}
As usual, $C$, $c$ and all constants in the proof might depend on
$C_{1}$ and $c_{1}$.

\begin{proof}
Clearly we may assume $E=1$. Let $Y_{i}$ be i.i.d variables with
$Y_{i}\sim C_{2}(1+G)$ where $G$ is a standard exponential variable
(namely with density $e^{-t}$) and $C_{2}$ is some constant sufficiently
large such that $\mathbb{P}(Y_{i}>\lambda )\geq \min \{1,C_{1}e^{-c_{1}\lambda }\}$.
A simple induction now shows that \[
\mathbb{P}\left(\sum X_{i}>\lambda n\right)\leq \mathbb{P}\left(\sum Y_{i}>\lambda n\right)\]
and the sum of the $Y_{i}$ has the distribution $C_{2}(n+\Gamma )$
where $\Gamma $ has density $e^{-t}\frac{t^{n-1}}{(n-1)!}$. A simple
calculation shows that $\mathbb{P}(\Gamma >\lambda n)\leq Ce^{-c\lambda }$.
\end{proof}
\begin{lem}
\label{lem:Ndr2md}Let $v\in T_{N}^{d}$, and let $R$ be a random
walk starting from $v$ going to a length of $C_{1}N^{d}r^{2-d}$
for some $C_{1}$ sufficiently large. Then \[
\mathbb{P}(R\textrm{ hits }B(0,r))\geq \frac{1}{2}\quad \forall v\in T_{N}^{d}.\]

\end{lem}
\begin{proof}
Define $S$ to be the preimage of $B$ in $\mathbb{Z}^{d}$ (namely
$B+N\mathbb{Z}^{d}$) and let $R'$ be a random walk of on $\mathbb{Z}^{d}$
starting from some preimage of $v$. Then\[
\mathbb{P}(R\textrm{ hits }B(0,r))=\mathbb{P}(R'[0,N^{d}r^{2-d}]\cap S\neq \emptyset ).\]
Define stopping times $t_{i}$ as follows: $t_{0}=0$ and for every
$i$ let $z_{i}$ be the element of $N\mathbb{Z}^{d}$ closest to
$R(t_{i})$. Define inductively \[
t_{i+1}=\min \{t>t_{i}:R(t)\in \partial B(z_{i},r)\cup \partial B(z_{i},2N)\}.\]
Since $d(R(t_{i}),z_{i})\leq N$ then using the harmonic potential
at $z_{i}$ shows that \[
\mathbb{P}(R(t_{i+1})\in S)\geq \mathbb{P}(R(t_{i+1})\in \partial B(z_{i},r)\geq c(r/N)^{d-2}\]
independently of the value of $R(t_{i})$. This immediately gives
that, for $C_{2}$ sufficiently large \begin{equation}
\mathbb{P}(R[0,t_{C_{2}(N/r)^{d-2}}]\cap S=\emptyset )\leq \left(1-c\left(\frac{r}{N}\right)^{d-2}\right)^{C_{2}(N/r)^{d-2}}\leq \frac{1}{4}.\label{eq:intersectquarter}\end{equation}
On the other hand, it is easy to see that $\mathbb{P}(t_{i+1}-t_{i}>\lambda N^{2})\leq Ce^{-c\lambda }$
independently of the past, and using lemma \ref{lem:sumexpexp} we
get that $\mathbb{P}(t_{n}>\lambda nN^{2})\leq Ce^{-c\lambda }$.
Using this for $n=C_{2}(N/r)^{d-2}$ and $\lambda $ sufficiently
large we get that \begin{equation}
\mathbb{P}(t_{n}>CN^{d}r^{2-d})\leq \frac{1}{4}.\label{eq:toobigquarter}\end{equation}
(\ref{eq:intersectquarter}) and (\ref{eq:toobigquarter}) together
show that the $C$ in (\ref{eq:toobigquarter}) may serve as our $C_{1}$.
\end{proof}
\begin{lem}
\label{lem:condexpout}Let $v\in \partial B(0,2r)$ and $w\in \partial B(0,r)$.
Let $R$ be a random walk started from $v$ and conditioned to hit
$B(0,r)$ in $w$. Let $t$ be the hitting time. Then\[
\mathbb{P}(t>\lambda N^{d}r^{2-d})\leq Ce^{-c\lambda }\quad .\]

\end{lem}
The proof is identical to that of lemma \ref{lem:condexpin} with
the use of lemmas \ref{lem:hitp} and \ref{lem:erw} replaced by the
comments following them, respectively, and using lemma \ref{lem:Ndr2md}
to show that the probability to not hit a ball of radius $r$ after
$\lambda N^{d}r^{2-d}$ steps is $\leq Ce^{-c\lambda }$ and hence
the equivalent of lemma \ref{lem:rrw_smooth}. We omit the details.

\end{document}